\documentclass[12pt]{article}
\usepackage{amsmath, amsfonts}
\usepackage{tikz}
\numberwithin{equation}{section}
\newtheorem{theorem}{Theorem}[section]
\newtheorem{lemma}[theorem]{Lemma}

\newtheorem{definition}[theorem]{Definition}
\newtheorem{remark}[theorem]{Remark}

\newcommand{\qed}{\nolinebreak\hfill\vbox{\hrule\hbox{\vrule\kern3pt\vbox{\kern6pt}\kern3pt\vrule}\hrule}}
\newenvironment{pf}{{\noindent\bf Proof.}}{\qed\newline}

\newcommand{\lap}{\Delta}

\begin{document}

\title{Symmetry and Approximate Symmetry of a Nonlinear Elliptic Problem over a Ring }
\author{Alaa Akram Haj Ali
\footnote{Alaa Akram Haj Ali is partially
supported by a GAAN Grant from the Department of Education.}\\
\footnotesize Department of Mathematics\\
\footnotesize Wayne State University\\
\footnotesize 656 W.Kirby, 1150 FAB\\
\footnotesize Detroit, MI 48202
\and Dongsheng Li
\footnote{Dongsheng Li is partially supported by NSFC 11671316}\\
\footnotesize Department of Mathematics\\
\footnotesize Xi'an Jiaotong University\\
\footnotesize 28 W. Xianning Road\\
\footnotesize Xi'an, China 710049
\and \& Peiyong Wang
\footnote{Peiyong Wang is partially supported by a Simon's Collaboration Grant.}\\
\footnotesize Department of Mathematics\\
\footnotesize Wayne State University\\
\footnotesize 656 W.Kirby, 1150 FAB\\
\footnotesize Detroit, MI 48202\\
\footnotesize E-mail: pywang@math.wayne.edu
\normalsize}
\date{}
\maketitle
\begin{abstract}
A singularly perturbed free boundary problem arising from a real problem associated with a Radiographic Integrated Test Stand concerns a solution of the equation $\Delta u = f(u)$ in a domain $\Omega$ subject to constant boundary data, where the function $f$ in general is not monotone. When the domain $\Omega$ is a perfect ring, we incorporate a new idea of radial correction into the classical moving plane method to prove the radial symmetry of a solution. When the domain is slightly shifted from a ring, we establish the stability of the solution by showing the approximate radial symmetry of the free boundary and the solution. For this purpose, we complete the proof via an evolutionary point of view, as an elliptic comparison principle is false, nevertheless a parabolic one holds.
\end{abstract}



\section{Introduction}\label{introduction}

Let $\Omega$ be the domain between two concentric spheres $|x| = 1$ and $|x| = R$ for some large radius $R$. Assume $u\in C^2(\bar{\Omega})$ is a solution of the boundary value problem
\begin{equation}\label{bvp}
\left\{\begin{array}{ll} \lap u = f(u) &\ \ \text{in $\Omega = B_R\backslash \bar{B}_1$,}\\
u = 1 &\ \ \text{on $|x| = 1$,} \\
u = -1 &\ \ \text{on $|x| = R$.}\end{array}\right.
\end{equation}
The function $f\colon\mathbb{R}^+\rightarrow\mathbb{R}$ is a $C^1$ function satisfying $f(s) \leq 0$. We study the radial symmetry of a solution of this boundary value problem. This problem is a singularly perturbed problem of the following free boundary problem arising from industry (cf. \cite{O}) the study of which will be the content of another paper:
\begin{equation}\label{fbp}
\left\{\begin{array}{ll}
\lap u = f(u)\ \ &\text{in\ } \left\{u > 0\right\} \\
\lap u = 0\ \ &\text{in\ } \left\{u \leq 0\right\} \\
u^+_{\nu} = u^-_{\nu}\ \ &\text{along\ } \mathcal{F}:= \partial\left\{u>0\right\}\\
u=1\ \ &\text{on\ } |x|\leq 1\\
u=-1\ \ &\text{on\ } |x| = R \end{array}\right.
\end{equation}
where $f(s) < 0$ for $s > 0$. If one allows $f(s) = 0$ when $s\leq 0$ and uses a smooth function to approximate this new function $f$, then one ends up with the problem (\ref{bvp}). We also consider the problem (\ref{bvp}) when the bigger sphere shifts its center a little from the origin
\begin{equation}\label{shift}
\left\{\begin{array}{ll} \lap u = f(u) &\ \ \text{in $\Omega = B_R(Z)\backslash \bar{B}_1$,}\\
u = 1 &\ \ \text{on $|x| = 1$,} \\
u = -1 &\ \ \text{on $|x-Z| = R$,}\end{array}\right.
\end{equation}
where $|Z| = \delta$ is small. The boundary of the positive set $\mathcal{F} := \partial\left\{u>0\right\}$ in each problem is the \textit{free boundary} of a solution $u$.

The goal of this paper is to prove the radial symmetry of a solution of the boundary value problem (\ref{bvp}), and the approximate radial symmetry of the free boundary of a solution of problem (\ref{shift}), under a not-too-negative condition on $f'$. With regard to the first task, our situation is different from known results in that there is no uniqueness of a solution for the Dirichlet problem, which can easily be seen. For example, suppose $\lambda$ is an eigenvalue of $(-\lap)$ with an eigenfunction $w$ on the region $\Omega = B_R\backslash\bar{B}_r$. That is
\begin{equation*}
\left\{\begin{array}{clcl}\lap w &= &-\lambda w &\ \ \text{in $\Omega$}\\ w &= &0 &\ \ \text{on $\partial B_r\cup\partial B_R$}\end{array}\right.
\end{equation*}
If $u$ is a solution of the Dirichlet problem
\begin{equation*}
\left\{\begin{array}{clcl}\lap u &= &-\lambda u &\ \ \text{in $\Omega$}\\ u &= &1 &\ \ \text{on $\partial B_r$}\\ u &= &-1 &\ \ \text{on $\partial B_R$,}\end{array}\right.
\end{equation*}
so is $u+w$. This does not happen for the primary eigenvalue but occurs for other eigenvalues according to the classical Courant's nodal set theorem. Existing results of symmetry or approximate symmetry of a solution over a ring-like domain depends on the assumption that the right-hand-side $f$ is non-decreasing. The reader may refer to \cite{HPP}, \cite{HP} and the references therein. One consequence of the non-uniqueness of a solution to the Dirichlet problem is the absence of a comparison principle for the equation, which is remedied in the standard moving plane method by the radial monotonicity of a solution. However, absence of the monotonicity of the right-hand-side $f$ in the current situation puts the radial monotonicity of a solution over a ring in question. The standard moving plane method does not work until this issue is resolved. In this sense, our method as well as results are new in the study of radial symmetry of a solution and may be applied in a broader scope in studying symmetry problems.

The second task of securing approximate radial symmetry when the domain is shifted somewhat from a ring has practical meaning in that it causes technical disaster and shutdown of the system if the free boundary touches the interior sphere in a Radiographic Integrated Test Stand or RITS (\cite{O}), and this is possible as in practice the interior and exterior spheres can never be perfectly concentric especially when RITS is in operation. In order to prove the approximate symmetry of a solution when the domain is shifted from a ring, we are, in a sense, forced to employ a technique of using evolutionary limits to bound the solution. The reason is the lack of an elliptic comparison principle and the uniqueness of a solution as stated above, and meanwhile we come to realize the validity of a parabolic comparison principle. We have not seen such an approach in the literature except the joint work \cite{CW} of one of the authors with Luis A.\,Caffarelli, in which the authors use a similar evolutionary view to examine the stability of a solution of an elliptic free boundary problem. Construction of the evolutionary limits depends on an existence theorem of a solution for the corresponding parabolic initial-boundary-value problem and locally uniform convergence of the evolution. In proving the existence theorem for an evolution, we are helped with an iteration rather than the widely used Perron's method, since the solution produced from that method may not be regular enough. This evolutionary approach to a problem in a steady state seems promising to us in application in the study of other PDE or free boundary problems.

The main results of this paper are the following two theorems regarding to problem (\ref{bvp}) and (\ref{shift}).
\begin{theorem}\label{symmetry}
Let $R > 1$ and $\Omega = B_R\backslash \bar{B}_1$ be the domain of a ring or shell. Suppose $f\colon \mathbb{R}_+\rightarrow\mathbb{R}$ is a $C^1$ function such that
$f(s)\leq 0$ and $\inf_{\mathbb{R}_+}f'(s) > -\frac{4(n+2)}{R^2}$.

Then a solution $u\in C^2(\overline{B_R\backslash B_1})$ of (\ref{bvp}) is radially symmetric in the sense $u(x) = u(y)$ if $x$, $y\in\Omega$ with $|x| = |y|$.
\end{theorem}

The definition of a \textit{stable solution} in the statement of the second theorem is given in Definition \ref{stable}.
\begin{theorem}\label{appr_symmetry}
Suppose $R > 1$, and $f\colon \mathbb{R}_+\rightarrow\mathbb{R}$ is a $C^1$ function such that $f(s)\leq 0$ and $\inf_{\mathbb{R}_+}f'(s) > -\frac{2(n+2)}{R^2}$. Let $u\in C^2(\overline{B_R(Z)\backslash B_1})$ be a stable solutions of (\ref{shift}) with free boundary $\mathcal{F}$, where $|Z| = \delta$.

Then there exists a constant $\delta_0 > 0$ such that for every constant $\delta$ in $0 < \delta \leq \delta_0$, there is a solution $u_0\in C^2(\overline{B_R\backslash\bar{B}_1})$ of (\ref{bvp}) with free boundary $\mathcal{F}_0$ so that
$$|u(x) - u_0(x)| \leq C\delta\ \ \text{in $\left(B_R(Z)\cap B_R\right)\backslash B_1$, and}$$
$$dist(\mathcal{F}, \mathcal{F}_0) < C|Z| = C\delta$$
for a constant $C = C(n,R,\inf f')$ which is independent of $\delta$. The latter estimate, in other words, states that the free boundary $\mathcal{F}$ is in the shell between two concentric spheres of thickness $2C\delta$, as Theorem \ref{symmetry} implies $\mathcal{F}_0$ is a sphere. In particular, the free boundary $\mathcal{F}$ keeps a positive distance from the boundary of the domain $\partial\Omega$.
\end{theorem}

In accordance with the goals, the rest of the paper is naturally divided into two parts. The next is devoted to the proof of the radial symmetry of a solution over a ring by way of the moving plane method. The third part presents the approximate symmetry when the domain is a shifted ring, in which well-posedness of the parallel evolution, convergence of the evolution, and bounds by the evolutionary limit solutions are established.

\section{Symmetry over a Ring}
In this section, one considers the following boundary value problem.
\begin{equation}\label{singularly_perturbed}
\left\{\begin{array}{cl} \lap u = f(u) &\ \text{in $1\leq |x| \leq R$}\\ u = 1 &\ \text{on $|x| = 1$}\\ u = 0 &\ \text{on $|x| = R$}\end{array}\right.
\end{equation}
One assumes $R$ is large, $u\in C^2(\overline{B_R\backslash B_1})$, and $f\colon \mathbb{R}_+\rightarrow\mathbb{R}$ is a $C^1$ function such that
$f(s)\leq 0$ and $\inf_{\mathbb{R}_+}f'(s) > -\frac{2(n+2)}{R^2}$. Let $\Omega = B_R\backslash \bar{B}_1$ be the domain of a ring or shell. We note the non-essential difference in the boundary value of a solution between the problems \ref{bvp} and \ref{singularly_perturbed}.

The goal of this section is to prove Theorem \ref{symmetry} which is equivalent to the following theorem.
\begin{theorem}\label{thm_perturbed}
Let $R > 1$ and $\Omega = B_R\backslash \bar{B}_1$ be the domain of a ring or shell. Suppose $f\colon \mathbb{R}_+\rightarrow\mathbb{R}$ is a $C^1$ function such that
$f(s)\leq 0$ and $\inf_{\mathbb{R}_+}f'(s) > -\frac{2(n+2)}{R^2}$.

Then a solution $u\in C^2(\overline{B_R\backslash B_1})$ of (\ref{singularly_perturbed}) is radially symmetric in the sense $u(x) = u(y)$ if $x$, $y\in\Omega$ with $|x| = |y|$.
\end{theorem}

\begin{remark}
The standard moving-plane argument, e.\,g.\,\cite{GNN}, stops in the middle sphere of the ring and hence cannot reach the radial symmetry. Moreover, as indicated in the Introduction, it is unknown if $u$ enjoys radial monotonicity. So a direct application of the moving-plane method does not work.
\end{remark}

We want to caution the reader that in general $u$ is not necessarily radially monotone. It could happen that $u$ assumes its maximum on a sphere in the ring $\Omega$ while still staying radially symmetric. This issue helps one to understand why we need to play the trick of adding a dominating radially symmetric function to enforce a radial symmetry on the resulting sum function, which we will describe now.

Firstly, one constructs an auxiliary dominating radially symmetric function. For any number $A > 0$, an alternating sequence $\left\{a_k\right\}^{\infty}_{k=0}$ is defined recursively by
\begin{equation*}
a_0 > 0,\ \ a_{k+1} = -\frac{Aa_k}{2(n+2k)(k+1)}.
\end{equation*}
One defines an analytic function $\phi$ on $\mathbb{R}$ by a power series
\begin{equation*}
\phi(s) = \sum^{\infty}_{k=0}a_ks^{2k},
\end{equation*}
which is obviously uniformly convergent on any bounded subset of $\mathbb{R}$. A direct computation shows that
\begin{equation*}
\phi''(s) + \frac{n-1}{s}\phi'(s) = -A\phi(s)\ \ (s\in\mathbb{R}),
\end{equation*}
which implies
\begin{equation*}
\lap \phi(|x|) = -A\phi(|x|)\ \ (x\in\mathbb{R}^n\backslash\left\{0\right\}).
\end{equation*}
In addition,
\begin{equation*}
\begin{split}
\phi'(s) &= \sum^{\infty}_{j=1}2(2j-1)a_{2j-1}\left(1 - \frac{As^2}{2(n+4j-2)(2j-1)}\right) s^{4j-3}\\
&< 0\ \ \ \ \text{if\ }\ s<\sqrt{\frac{2(n+2)}{A}}
\end{split}
\end{equation*}
Moreover, if one requires $$-\inf_{\mathbb{R}_+}f'(s) < A < \frac{2(n+2)}{R^2},$$ then for $s \leq R$ it holds
\begin{equation*}
\begin{split}
\phi'(s) &\leq 2a_1\left(1 - \frac{As^2}{2(n+2)}\right)\\
&\leq -\frac{Aa_0}{n}\left(1 - \frac{AR^2}{2(n+2)}\right)
\end{split}
\end{equation*}

We will apply the well-known moving plane method which plays the key role in \cite{S} and \cite{GNN} to the function
\begin{equation}\label{sum_function}
\tilde{u}(x) = u(x) + C\phi(|x|)
\end{equation}
in $\Omega$ for positive constants $A$ and $C$. We pick the value of $C$ so that $$C\geq \frac{n}{Aa_0\left(1-\frac{AR^2}{2(n+2)}\right)} \sup_{\Omega}\left| \nabla u(x)\right|.$$ Then $\tilde{u}_r(x) \leq 0$ for all $x\in \Omega$, i.\,e. $\tilde{u}$ is radially decreasing.

For any domain $\mathcal{D}$ in consideration, $\nu(x_0)$ denotes the outer unit normal to $\partial\mathcal{D}$ at a point $x_0\in\partial\mathcal{D}$.

In order to prove $u$ is radially symmetric in $\Omega$, it suffices to prove $\tilde{u}$ is radially symmetric in the ring $\Omega$, which is equivalent to that $\tilde{u}$ is symmetric in every hyperplane through the origin. Without loss of generality, one takes the direction $\nu = e_1$ and starts to prove $\tilde{u}$ is symmetric in the hyperplane $x_1 = 0$.

For the sake of completeness of this work, we include here the version of Hopf's lemma and Strong Maximum Principle that we will use in the proof.

\begin{theorem}
\textbf{Hopf's Lemma}

Suppose $u\in C^2(\Omega)\cap C^1(\bar{\Omega})$ is a solution of the differential inequality $$\lap u(x) + c(x)u(x) \geq 0$$ in $\Omega$, where $c\in C(\Omega)$. Assume further $u(x) < 0$ in $\Omega$, $x_0\in\partial\Omega$ such that $u(x_0) = 0$, and there is a ball $B\subset\Omega$ that touches $\partial\Omega$ at $x_0$.

Then $$u_{\nu}(x_0) > 0$$ for the unit outer normal $\nu$ at $x_0$ to $\partial\Omega$.
\end{theorem}
For a proof of the Hopf's lemma, the reader may refer to \cite{E} for the case $c(x)\leq 0$ and \cite{GNN} for the case $c(x) > 0$.

\begin{theorem}
\textbf{Strong Maximum Principle}

Suppose $\Omega$ is connected and $u\in C^2(\Omega)$ is a solution of the differential inequality $$\lap u(x) + c(x)u(x) \geq 0$$ in $\Omega$, where $c\in C(\Omega)$, and $u(x)\leq 0$ in $\Omega$.

If $u(x_0) = 0$ at a point $x_0$ in $\Omega$, then $u(x) \equiv 0$ in $\Omega$.
\end{theorem}

\begin{figure}
\begin{tikzpicture}
<\fill[blue!70!white] (1.5,-1.3229) arc (-41.41:41.41:2cm); \filldraw[step=0.25cm, gray!25!white, ultra thin] (-2.47,-2.47) grid (2.47,2.47); \draw[step=1cm, gray!75!white, very thin] (-2.47,-2.47) grid (2.47,2.47); \draw[thick,black, ->] (-2.47,0) -- (2.47,0) node[anchor=north] {\tiny x}; \draw[thick,black, ->] (0,-2.47) -- (0, 2.47) node[anchor=east] {\tiny y};
\node[anchor=north] at (-2,0) {\tiny -2}; \node[anchor=north] at (-1,0) {\tiny -1}; \node[anchor=north east] at (0,0) {\tiny 0}; \node[anchor=north] at (1,0) {\tiny 1}; \node[anchor=north] at (2,0) {\tiny 2};
\draw[black] (0,0) circle (1cm); \draw[black] (0,0) circle (2cm);
\node[anchor=east] at (0,-2) {\tiny -2}; \node[anchor=east] at (0,-1) {\tiny -1}; \node[anchor=east] at (0,1) {\tiny 1}; \node[anchor=east] at (0,2) {\tiny 2};
\draw[black] (0,0) circle (1cm); \draw[black] (0,0) circle (2cm);>
\end{tikzpicture}
\begin{tikzpicture}
<\fill[blue!70!white] (1.25,-1.5612) arc (-51.32:51.32:2cm); \fill[white] (2,-0.8660) arc (240:120:1cm); \filldraw[step=0.25cm, gray!25!white, ultra thin] (-2.47,-2.47) grid (2.47,2.47); \draw[step=1cm, gray!75!white, very thin] (-2.47,-2.47) grid (2.47,2.47); \draw[thick,black, ->] (-2.47,0) -- (2.47,0) node[anchor=north] {\tiny x}; \draw[thick,black, ->] (0,-2.47) -- (0, 2.47) node[anchor=east] {\tiny y};
\node[anchor=north] at (-2,0) {\tiny -2}; \node[anchor=north] at (-1,0) {\tiny -1}; \node[anchor=north east] at (0,0) {\tiny 0}; \node[anchor=north] at (1,0) {\tiny 1}; \node[anchor=north] at (2,0) {\tiny 2};
\draw[black] (0,0) circle (1cm); \draw[black] (0,0) circle (2cm);
\node[anchor=east] at (0,-2) {\tiny -2}; \node[anchor=east] at (0,-1) {\tiny -1}; \node[anchor=east] at (0,1) {\tiny 1}; \node[anchor=east] at (0,2) {\tiny 2};
\draw[black] (0,0) circle (1cm); \draw[black] (0,0) circle (2cm);>
\end{tikzpicture}
\begin{tikzpicture}
<\fill[blue!70!white] (0.75,-1.8540) arc (-67.98:67.98:2cm); \fill[white] (1.5,0) circle (1cm); \filldraw[step=0.25cm, gray!25!white, ultra thin] (-2.47,-2.47) grid (2.47,2.47); \draw[step=1cm, gray!75!white, very thin] (-2.47,-2.47) grid (2.47,2.47); \draw[thick,black, ->] (-2.47,0) -- (2.47,0) node[anchor=north] {\tiny x}; \draw[thick,black, ->] (0,-2.47) -- (0, 2.47) node[anchor=east] {\tiny y};
\node[anchor=north] at (-2,0) {\tiny -2}; \node[anchor=north] at (-1,0) {\tiny -1}; \node[anchor=north east] at (0,0) {\tiny 0}; \node[anchor=north] at (1,0) {\tiny 1}; \node[anchor=north] at (2,0) {\tiny 2};
\draw[black] (0,0) circle (1cm); \draw[black] (0,0) circle (2cm);
\node[anchor=east] at (0,-2) {\tiny -2}; \node[anchor=east] at (0,-1) {\tiny -1}; \node[anchor=east] at (0,1) {\tiny 1}; \node[anchor=east] at (0,2) {\tiny 2};
\draw[black] (0,0) circle (1cm); \draw[black] (0,0) circle (2cm);>
\end{tikzpicture}
\begin{tikzpicture}
<\fill[blue!70!white] (0.5,-1.9365) arc (-75.52:75.52:2cm); \fill[white] (1,0) circle (1cm); \filldraw[step=0.25cm, gray!25!white, ultra thin] (-2.47,-2.47) grid (2.47,2.47); \draw[step=1cm, gray!75!white, very thin] (-2.47,-2.47) grid (2.47,2.47); \draw[thick,black, ->] (-2.47,0) -- (2.47,0) node[anchor=north] {\tiny x}; \draw[thick,black, ->] (0,-2.47) -- (0, 2.47) node[anchor=east] {\tiny y};
\node[anchor=north] at (-2,0) {\tiny -2}; \node[anchor=north] at (-1,0) {\tiny -1}; \node[anchor=north east] at (0,0) {\tiny 0}; \node[anchor=north] at (1,0) {\tiny 1}; \node[anchor=north] at (2,0) {\tiny 2};
\draw[black] (0,0) circle (1cm); \draw[black] (0,0) circle (2cm);
\node[anchor=east] at (0,-2) {\tiny -2}; \node[anchor=east] at (0,-1) {\tiny -1}; \node[anchor=east] at (0,1) {\tiny 1}; \node[anchor=east] at (0,2) {\tiny 2};
\draw[black] (0,0) circle (1cm); \draw[black] (0,0) circle (2cm);>
\end{tikzpicture}
\begin{tikzpicture}
<\fill[blue!70!white] (0.25,-1.9843) arc (-82.82:82.82:2cm); \fill[white] (0.5,0) circle (1cm); \filldraw[step=0.25cm, gray!25!white, ultra thin] (-2.47,-2.47) grid (2.47,2.47); \draw[step=1cm, gray!75!white, very thin] (-2.47,-2.47) grid (2.47,2.47); \draw[thick,black, ->] (-2.47,0) -- (2.47,0) node[anchor=north] {\tiny x}; \draw[thick,black, ->] (0,-2.47) -- (0, 2.47) node[anchor=east] {\tiny y};
\node[anchor=north] at (-2,0) {\tiny -2}; \node[anchor=north] at (-1,0) {\tiny -1}; \node[anchor=north east] at (0,0) {\tiny 0}; \node[anchor=north] at (1,0) {\tiny 1}; \node[anchor=north] at (2,0) {\tiny 2};
\draw[black] (0,0) circle (1cm); \draw[black] (0,0) circle (2cm);
\node[anchor=east] at (0,-2) {\tiny -2}; \node[anchor=east] at (0,-1) {\tiny -1}; \node[anchor=east] at (0,1) {\tiny 1}; \node[anchor=east] at (0,2) {\tiny 2};
\draw[black] (0,0) circle (1cm); \draw[black] (0,0) circle (2cm);>
\end{tikzpicture}
\begin{tikzpicture}
<\fill[blue!70!white] (0,-2) arc (-90:90:2cm); \fill[white] (0,0) circle (1cm); \filldraw[step=0.25cm, gray!25!white, ultra thin] (-2.47,-2.47) grid (2.47,2.47); \draw[step=1cm, gray!75!white, very thin] (-2.47,-2.47) grid (2.47,2.47); \draw[thick,black, ->] (-2.47,0) -- (2.47,0) node[anchor=north] {\tiny x}; \draw[thick,black, ->] (0,-2.47) -- (0, 2.47) node[anchor=east] {\tiny y};
\node[anchor=north] at (-2,0) {\tiny -2}; \node[anchor=north] at (-1,0) {\tiny -1}; \node[anchor=north east] at (0,0) {\tiny 0}; \node[anchor=north] at (1,0) {\tiny 1}; \node[anchor=north] at (2,0) {\tiny 2};
\draw[black] (0,0) circle (1cm); \draw[black] (0,0) circle (2cm);
\node[anchor=east] at (0,-2) {\tiny -2}; \node[anchor=east] at (0,-1) {\tiny -1}; \node[anchor=east] at (0,1) {\tiny 1}; \node[anchor=east] at (0,2) {\tiny 2};
\draw[black] (0,0) circle (1cm); \draw[black] (0,0) circle (2cm);>
\end{tikzpicture}
\caption{$\Pi(\lambda)$ for $R = 2$, $\lambda = 1.5, 1.25, 1, 0.75, 0.5, 0.25, 0$, respectively}\label{fig1}
\end{figure}
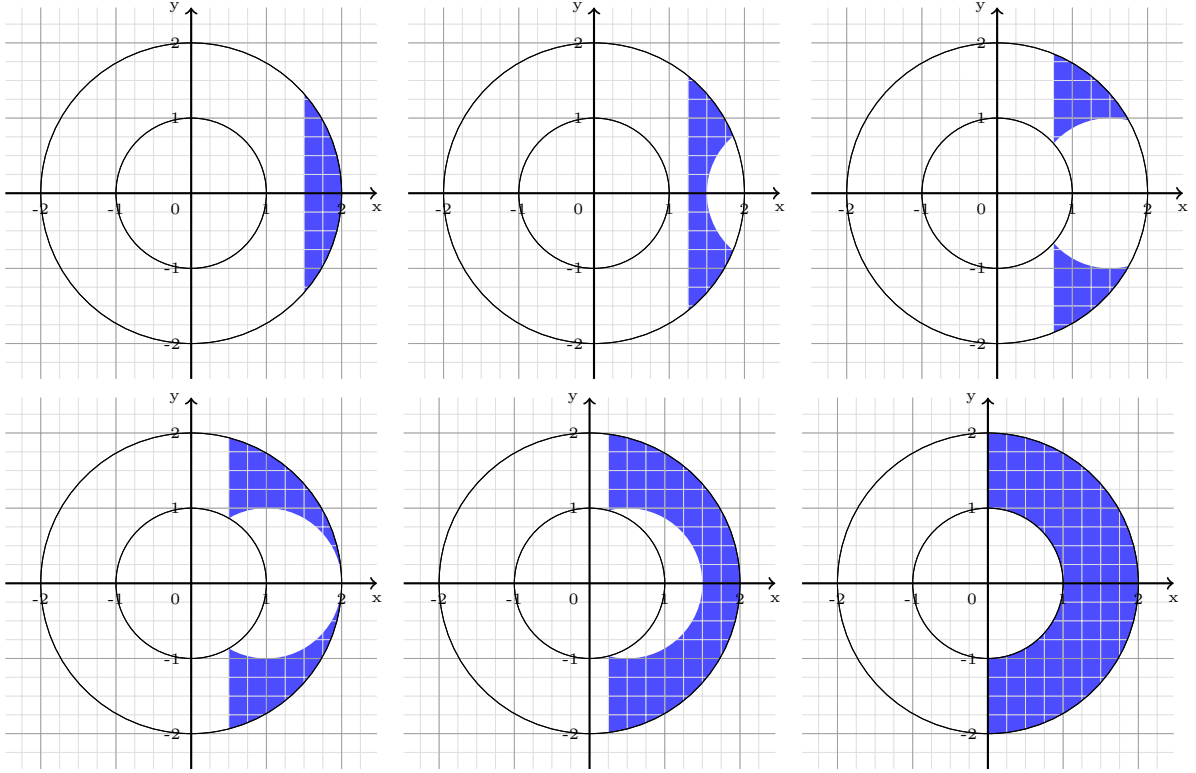

For any $\lambda \geq 0$, let $T_{\lambda}$ be the hyperplane $x_1 = \lambda$, $x^{\lambda} = (2\lambda - x_1, x_2, \hdots, x_n)$ be the mirror image of $x = (x_1, x_2, \hdots, x_n)$ in $T_{\lambda}$, $\Sigma(\lambda) = \Omega\cap \left\{x\colon x_1 > \lambda\right\}$, $\Pi(\lambda) = \left\{x\in\Sigma(\lambda)\colon x^{\lambda}\in\Omega\right\}$, $\Sigma'(\lambda)$ the reflection of $\Sigma(\lambda)$ in $T_{\lambda}$, and $\Pi'(\lambda) = \Sigma'(\lambda)\cap\Omega$ the reflection of $\Pi(\lambda)$ in $T_{\lambda}$. Figure \ref{fig1} provides some snapshots of the domain $\Pi(\lambda)$, shaded in blue, during the motion of the hyperplane at different values of $\lambda$ when the outer radius of the ring $R = 2$.

If one notices that $u$ is super-harmonic in $\Omega$ and attains its minimum on the sphere $|x| = R$, it is obvious the following lemma is true.
\begin{lemma}\label{lemma1}
Suppose $x_0\in\partial B_R$ with $\nu_1(x_0) > 0$.

Then there exists $\delta > 0$ such that
$$u_{x_1} < 0\ \ \text{and hence\ }\ \tilde{u}_{x_1} < 0$$ in $\Omega\cap\left\{x\colon |x-x_0| < \delta\right\}$.
\end{lemma}

The next lemma allows one to move the hyperplane $T_{\lambda}$ for $\lambda > 0$ in the negative $x_1$-axis direction.
\begin{lemma}\label{lemma2}
Fix some $\lambda$ in $0\leq \lambda < R$. Assume
$$\tilde{u}_{x_1}(x)\leq 0\ \ \ \text{in\ } \Sigma(\lambda)\ \ \ \text{and\ } \tilde{u}(x) \leq \tilde{u}(x^{\lambda})\ \ \text{in\ } \Pi(\lambda),$$ but $\tilde{u}(x)\not\equiv \tilde{u}(x^{\lambda})$ in $\Pi(\lambda)$.

Then $\tilde{u}(x) < \tilde{u}(x^{\lambda})$ in $\Pi(\lambda)$ and $\tilde{u}_{x_1}(x) < 0$ on $\Omega\cap T_{\lambda}$.
\end{lemma}
\begin{pf}
On $\overline{\Pi'(\lambda)}$, one defines the functions
\begin{equation*}
\begin{split}
&v(x) = u(x^{\lambda}),\ \
\tilde{v}(x) = \tilde{u}(x^{\lambda}) = u(x^{\lambda}) + C\phi(|x^{\lambda}|),\\
\text{and\ }\ &h(x) = C\phi(|x^{\lambda}|) - C\phi(|x|) \leq 0.
\end{split}
\end{equation*}

Define $w(x) = \tilde{v}(x) - \tilde{u}(x)$ on $\overline{\Pi'(\lambda)}$. Then $w(x) \leq 0$ in $\Pi'(\lambda)$ and $w$ satisfies
\begin{equation*}
\lap w + c(x)w = -\int^1_0f'((1-t)u + tv)\,dt\,h + \lap h
\end{equation*}
for $$c(x) = -\int^1_0f'((1-t)u+tv)\,dt$$ which is a continuous function on $\Omega$,
due to the equality
\begin{equation*}
\begin{split}
\lap(v - u + h) &= f(v) - f(u) + \lap h \\
&= \int^1_0f'((1-t)u+tv)\,dt\,(v-u) + \lap h.
\end{split}
\end{equation*}

As a consequence,
\begin{equation*}
\begin{split}
\lap w + c(x)w &\geq -\inf_{\mathbb{R}}f'(s)\,h + \lap h \\
&\geq Ah + \lap h \\
&= 0
\end{split}
\end{equation*}
as
\begin{equation*}
\begin{split}
\lap h(x) &= \lap \left(C\phi(|x^{\lambda}|)\right) - \lap \left(C\phi(|x|)\right) = -AC\phi(|x^{\lambda}|) + AC\phi(|x|) \\
&= - Ah(x).
\end{split}
\end{equation*}
Notice that $w(x) = 0$ on $T_{\lambda}\cap\bar{\Omega}$ and $w(x) \leq 0$ elsewhere on $\partial\Pi'(\lambda)$. Then the Strong Maximum Principle implies $w < 0$ in $\Pi'(\lambda)$, and the Hopf's Lemma implies $w_{x_1}(x) > 0$ on $T_{\lambda}\cap\Omega$. These mean
$$\tilde{v}(x) < \tilde{u}(x)\ \ \text{in $\Pi'(\lambda)$, or equivalently\ } \tilde{u}(x) < \tilde{u}(x^{\lambda})\ \ \text{in $\Pi(\lambda)$}$$ and
$\tilde{u}_{x_1}(x) < 0$ on $\Omega\cap T_{\lambda}$, since $w_{x_1}(x) = -\tilde{u}_{x_1}(x^{\lambda}) - \tilde{u}_{x_1}(x) = -2\tilde{u}_{x_1}(x)$ on $T_{\lambda}\cap\Omega$.
\end{pf}

The main Theorem (\ref{thm_perturbed}) follows from the following theorem by considering all possible directions along which a hyperplane is moved.
\begin{theorem}\label{thmlap}
For any $\lambda$ in $0 < \lambda < R$, it holds that
\begin{equation}\label{transient-condition}
\tilde{u}_{x_1}(x) < 0 \ \ \text{in $\Sigma(\lambda)$}\ \ \text{and\ } \tilde{u}(x) < \tilde{u}(x^{\lambda})\ \ \text{in $\Pi(\lambda)$.}
\end{equation}
In particular, $\tilde{u}_{x_1}(x) < 0$ in $\Omega\cap\left\{x_1>0\right\}$.

Consequently, $\tilde{u}(x)$ is symmetric with respect to the hyperplane $x_1 =0$.
\end{theorem}
\begin{pf}
We define the set $\mathcal{A}$ as
\begin{equation*}
\mathcal{A} = \left\{\lambda\in (0, R) \colon \tilde{u}_{x_1}(x) < 0\ \ \text{in $\Sigma(\lambda)$}\ \ \text{and\ \ } \tilde{u}(x) < \tilde{u}(x^{\lambda})\ \ \text{in $\Pi(\lambda)$} \right\}.
\end{equation*}
Firstly, one notices that Lemma \ref{lemma1} implies there exists some $\lambda$ close to $R$ in $0 < \lambda < R$ which is in $\mathcal{A}$.

Let $\mu = \inf\mathcal{A}$. Since (\ref{transient-condition}) holds for all $\lambda > \mu$, we have by continuity that
$$\tilde{u}_{x_1}(x) < 0\ \ \text{in $\Sigma(\mu)$}\ \ \text{and\ \ } \tilde{u}(x)\leq \tilde{u}(x^{\lambda})\ \ \text{in $\Pi(\mu)$.}$$
We claim that $\mu = 0$.

Suppose $\mu > 0$. For any $x_0\in\left(\partial B_R\cap\left\{x_1>\mu\right\}\right)$ such that $x^{\mu}_0\in\Omega$, it holds that $-1 + \phi(R) = \min_{\bar{\Omega}}\tilde{u} =  \tilde{u}(x_0) < \tilde{u}(x^{\mu}_0)$. So $\tilde{u}(x)\not\equiv \tilde{u}(x^{\lambda})$ in $\Pi(\mu)$. Lemma \ref{lemma2} then implies
$$\tilde{u}(x) < \tilde{u}(x^{\mu})\ \ \text{in\ } \Pi(\mu)\ \ \text{and\ } \tilde{u}_{x_1}(x) < 0\ \ \text{on\ } \Omega\cap T_{\mu}.$$ That is, (\ref{transient-condition}) holds for $\lambda = \mu$.

At every point $x_0\in\partial\Omega\cap T_{\mu}$, Lemma \ref{lemma1} states there is a $\varepsilon > 0$ such that
$$\tilde{u}_{x_1} < 0\ \ \text{in\ } \Omega\cap \left\{|x-x_0| < \varepsilon\right\},$$ as $T_{\mu}$ is not perpendicular to $\partial \Omega$. Here one notices that the situation when $|x_0| = 1$ is parallel to that in Lemma \ref{lemma1} and a similar conclusion holds. Since $\partial\Omega\cap T_{\mu}$ is compact, there is an $\varepsilon > 0$ such that
$$\tilde{u}_{x_1} < 0\ \ \text{in\ } \Omega\cap\left\{x_1 > \mu - \varepsilon\right\}\cap N_{\varepsilon}(\partial\Omega\cap T_{\mu}),$$ where $N_{\varepsilon}(S)$ denotes the $\varepsilon$-neighborhood of a set $S\in\mathbb{R}^n$.
On the other hand, since $\tilde{u}_{x_1} < 0$ on $\Omega\cap T_{\mu}$, one gets by continuity of $\tilde{u}_{x_1}$ that
$$\tilde{u}_{x_1} < 0\ \ \text{in\ } \Omega\cap \left\{x_1 > \mu - \varepsilon\right\}\backslash N_{\varepsilon}(\partial\Omega\cap T_{\mu})$$
so long as the value of $\varepsilon$ is taken smaller if necessary. In all, for this $\varepsilon > 0$,
\begin{equation}\label{2.5}
\tilde{u}_{x_1} < 0\ \ \text{in\ } \Omega\cap \left\{x_1 > \mu - \varepsilon\right\}.
\end{equation}

As $\mu = \inf\mathcal{A}$, $\exists\left\{\lambda^j\right\}$ such that $0 < \lambda^j < \mu$ and
$$\exists x_j\in \Pi(\lambda^j)\ \ \text{such that\ } \tilde{u}(x_j)\geq \tilde{u}(x^{\lambda^j}_j)\ \ \text{for every $j$.}$$
Without loss of generality, we assume $x_j\rightarrow \tilde{x}$ for some $\tilde{x}\in \overline{\Pi(\mu)}$. Clearly $x^{\lambda^j}_j\rightarrow \tilde{x}^{\mu}$ and hence $\tilde{u}(\tilde{x}) \geq \tilde{u}(\tilde{x}^{\mu})$. Since (\ref{transient-condition}) holds for $\lambda = \mu$, we must have $\tilde{x}\in \partial\Pi(\mu)$. There are four possibilities, $|\tilde{x}| = 1$, $|\tilde{x}| = R$, $\tilde{x}\in T_{\mu}\cap\Omega$, and $x\in \left(\partial \Pi(\mu)\backslash T_{\mu}\right)\cap\Omega$. One first notes that it is impossible that $|\tilde{x}| = 1$ but $\tilde{x}\not\in T_{\mu}$, since otherwise $\left|\left(x^j\right)^{\mu}\right| < 1$ holds for sufficiently large $j$ due to $\mu > 0$.
If $\left|\tilde{x}\right| = R$, then $\tilde{x}^{\mu}\in\Omega$ or $|\tilde{x}^{\mu}| = 1$, and since $\tilde{u}$ is radially decreasing,
$$\tilde{u}(\tilde{x}) = \min_{\bar{\Omega}}\tilde{u} < \tilde{u}(\tilde{x}^{\mu}),\ \ \text{which is a contradiction.}$$
Similarly, we get a contradiction when $\tilde{x}\in \left(\partial\Pi(\mu)\backslash T_{\mu}\right)\cap\Omega$, since, in this case, $|\tilde{x}^{\mu}| = 1$ and the fact $\tilde{u}$ is radially decreasing imply $$\tilde{u}(\tilde{x}) < \max_{\bar{\Omega}}\tilde{u} = \tilde{u}(\tilde{x}^{\mu}).$$
Therefore $\tilde{x}\in T_{\mu}\cap\bar{\Omega}$ and $\tilde{x}^{\mu} = \tilde{x}$. On the other hand, for large $j$, the segment $[x_j, x^{\lambda^j}_j]\subset\Omega$ and therefore $\exists y_j\in [x_j, x^{\lambda^j}_j]$ such that $u_{x_1}(y_j) \geq 0$ according to the Mean Value Theorem. Since $y_j\rightarrow \tilde{x}$, we get $u_{x_1}(\tilde{x}) \geq 0$ which is in contradiction to (\ref{2.5}).

Thus $\mu = 0$ and (\ref{transient-condition}) holds for all $\lambda$ in $0 < \lambda < R$. By continuity, it holds that $\tilde{u}_{x_1}(x) \leq 0$ and $\tilde{u}(x) \leq \tilde{u}(x^{0})$ in $\Sigma(0)$, where $x^0$ is the reflection of $x$ in the hyperplane $x_1 = 0$.

If one moves the hyperplane along the positive $x_1$-axis direction from the other side of the ring $\Omega$, the above argument shows that $\tilde{u}(x) \geq \tilde{u}(x^0)$ and hence $\tilde{u}$ and therefore $u$ are symmetric about the hyperplane $x_1 = 0$.
\end{pf}

The main theorem \ref{thm_perturbed} of this section follows readily from the preceding theorem.


\section{Stability of the Free Boundary}
This section is devoted to the proof of Theorem \ref{appr_symmetry}. Let $\Omega = B_R(Z)\backslash \bar{B}_1$ be a slight deformation of the ring $B_R\backslash \bar{B}_1$ with $|Z| = \delta > 0$ being sufficiently small. Now one considers the following boundary value problem.
\begin{equation}\label{singularly_perturbed_1}
\left\{\begin{array}{ll} \lap u = f(u) &\ \text{in $\Omega$}\\ u = 1 &\ \text{on $|x| = 1$}\\ u = -1 &\ \text{on $|x-Z| = R$}\end{array}\right.
\end{equation}
One assumes $R > 1$, $u\in C^2(\overline{\Omega})$, and $f\colon \mathbb{R}\rightarrow\mathbb{R}$ is a $C^3$ function such that $f(s)\leq 0$, $f(s)= 0$ if $s\leq 0$, and $\inf_{\mathbb{R}_+}f'(s) > -\frac{2(n+2)}{R^2}$. We consider only the stability of the free boundaries of what we call \textit{stable solutions} in a strong sense defined below.
\begin{definition}\label{stable}
A solution $u$ of (\ref{singularly_perturbed_1}) is \textbf{stable} if for any $\varepsilon > 0$, there exist functions $v_1$ and $v_2$ in $C^2(\bar{\Omega})$ that satisfy
\begin{equation}\label{3.2}
u - \varepsilon \leq v_1\leq u\leq v_2\leq u+ \varepsilon\ \ \text{on $\bar{\Omega}$,}
\end{equation}
\begin{equation}\label{3.3}
-\Delta v_1 + f(v_1) < -\varepsilon\ \ \text{and} \ \ -\Delta v_2 + f(v_2) > \varepsilon\ \ \text{in $\Omega$, simultaneously.}
\end{equation}
\end{definition}
\begin{remark}
When the domain is a ring and $f(s) \equiv 0$, it is easy to construct the sub- and super-solutions $v_1$ and $v_2$. One may readily perturb the domain to a ring-like one such as $\Omega$ and construct corresponding sub- and super-solutions over $\Omega$ that satisfy the requirements in the above definition. The reader is referred to the following proof for detailed computation.
\end{remark}

In other words, a stable solution $u$ is a uniform supremum of strict subsolutions and a uniform infimum of strict supersolutions. Compared to the concentric case when $Z = 0$, the center of the exterior sphere drifts away from the origin a bit. Our goal in this section is to prove in this situation the free boundary of $u$ drifts away from its original position also by a bit. In mathematical terms, we are to prove the stability of the free boundary. We will also give an estimate of the drift of the free boundary. 
However, for this seemingly clear fact, we need to prove it through a delicate evolution with quite a few technicalities. The reason we go through this quite troublesome process lies in the observation there is no comparison principle and hence no uniqueness for the elliptic problem when the nonlinear term $f(u)$ is negative. Nevertheless, there is a comparison principle for the corresponding evolution. Meanwhile, the reader may have realized that the practical reason why we study this problem on approximate radial symmetry has already been mentioned in the introduction.

We first state the parabolic comparison principle which is needed in the coming proof. Consider the initial-boundary value problem
\begin{equation}\label{ibv}
\left\{\begin{array}{ll}Hw := w_t - \Delta w + \alpha(x,w) = 0 &\ \text{in\ }\Omega\times (0, \infty)\\ w(x,t) = \sigma(x,t)\ \ \text{on\ }\partial\Omega\times (0,\infty),\ \ &\ w(x,0) = v_0(x)\ \ \text{for\ }x\in\bar{\Omega}\end{array}\right.
\end{equation}
where $\alpha$ is a $C^1$ function that satisfies the condition $0\leq \alpha(x,w)\leq Cw$, and $\Omega$ is a bounded domain with smooth boundary. This problem includes two important cases that we will apply the comparison principle to, the case when $\alpha = f(w)$ and the other when $\alpha = f'(w)z$ where $z$ is one of the first order derivatives of $w$.
\begin{theorem}\label{para_cp}

Suppose two functions $w_1$ and $w_2$ satisfy $Hw_1 \leq 0 \leq Hw_2$ in the viscosity sense as continuous functions or in the weak sense as $H^1$-functions in $\Omega\times\mathbb{R}^+$ and $w_1\leq w_2$ on the parabolic boundary $\partial_p(\Omega\times\mathbb{R}^+)$. Then $w_1\leq w_2$ in $\Omega\times\mathbb{R}^+$. Here $\mathbb{R}^+ = (0,\infty)$.
\end{theorem}
\begin{pf}
The proof is done with the introduction of the new functions $$\tilde{w}_j(x,t) = e^{-\lambda t}\left(w_j(x,t) - \frac{\delta}{T-t}\right),\ \ j= 1,2,$$ for any fixed small $T > 0$ and some large constant $\lambda$, cf.\,Theorem 3.1 \cite{CW} and Lemma 6.3 \cite{LW}.
\end{pf}

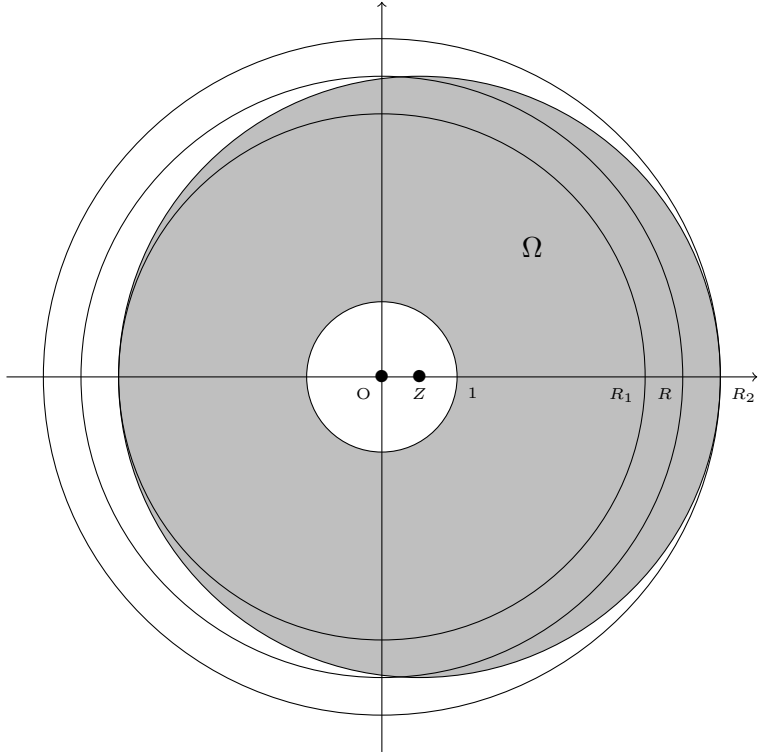
\begin{figure}
\begin{center}
\begin{tikzpicture}
<\fill[gray!50!white] (0.5,0) circle (4cm); \fill[white] (0,0) circle (1cm);
\draw[black] (0,0) circle (1cm); \draw[black] (0,0) circle (3.5cm); \draw[black] (0.5,0) circle (4cm); \draw[black] (0,0) circle (4cm); \draw[black] (0,0) circle (4.5cm);
\node[black] at (0,0) {\textbullet}; \node[anchor=north east] at (0,0) {\tiny O}; \node[black] at (0.5,0) {\textbullet}; \node[anchor=north] at (0.5,0) {\tiny $Z$};
\node[black,anchor=north east] at (3.5,0) {\tiny $R_1$}; \node[black,anchor=north east] at (4,0) {\tiny $R$}; \node[black,anchor=north west] at (4.5,0) {\tiny $R_2$}; \node[anchor=north] at (2,2) {\small $\Omega$}; \node[anchor=north west] at (1,0) {\tiny 1};
\draw[thin,black,->] (-4.99,0) -- (4.99,0); \draw[thin,black,->] (0,-4.99) -- (0, 4.99);
>
\end{tikzpicture}
\caption{The spheres $B_1$, $B_{R_1}$, $B_R(Z)$, $B_R$, and $B_{R_2}$ for $\delta = 0.5$ and $R = 4$}\label{fig2}
\end{center}
\end{figure}

Now let $B_{R_1}$ be the largest ball inscribed in $B_R(Z)$ with the origin as its center and $B_{R_2}$ be the smallest ball circumscribing $B_R(Z)$ with the origin as its center. Also, let $\mathcal{R} = B_R\backslash\bar{B}_1$ be a concentric ring, $\Omega_1=B_{R_1} \setminus \bar{B}_1$ and $\Omega_2=B_{R_2} \setminus \bar{B}_1$. Figure \ref{fig2} illustrates the two-dimensional sections of these spheres and the domain $\Omega$ as shaded in gray.

Let $u$ be a stable solution of the free boundary problem (\ref{singularly_perturbed_1}). Fix a small number $\varepsilon = K\delta$ for a relative large universal constant $K > 0$ in (\ref{3.2}) and (\ref{3.3}). Let $v_1$ and $v_2$ be as in the definition of the stable solution $u$ in $\Omega$. It is not difficult to see that, in accordance with the definition of $v_1$ and $v_2$, $v_1 < u < v_2$ on $\partial\Omega$.

In the following, we will construct a function $v_{01}$ (resp. $v_{00}$ and  $v_{02}$)  a strict subsolution (resp. strict supersolutions) of our problem on the perfect ring $\Omega_1 $ (resp.$\mathcal{R}$ and $\Omega_2 $) such that
\begin{equation}
\begin{split}
&u-C\delta \leq v_{01} \leq u  \text{ in } \Omega_1, \text{ and, } u \leq v_{02} \leq u+C\delta \text{ in } \Omega \\
&  v_{01} \leq v_{00} \text{ in }  \Omega_1\text { and, } v_{00} \leq v_{02} \text{ in }  \mathcal{R} \\
\end{split}
\end{equation}
for a constant $C$. 

Then we will use $v_{01}$ (resp.\,$v_{00}$ and $v_{02}$) as initial data of the parabolic version of our problem on $\Omega_1 \times (0,\infty)$ (resp. $\mathcal{R}\times (0,\infty)$ and $\Omega_2\times (0,\infty) $) to construct solutions of the respective evolution.

Finally, we prove convergence of the evolution with each initial data to a steady state which gives desired solutions $u_1$, $u_0$, and $u_2$ of the elliptic problems on $\Omega_1$, $\mathcal{R}$, and $\Omega_2$. The solutions $u_1$ and $u_2$ will give the lower and upper bounds for the solution $u$ of (\ref{shift}), while $u_0$ will be a radially symmetric approximation of $u$. In particular, the free boundary of $u_0$ is an approximation of that of $u$.

\subsection{Construction of a solution of our problem on the perfect ring $\Omega_1 $}
\subsubsection{Construction of a strict subharmonic function in $\Omega$ satisfying the boundary conditions associated with our problem}

One takes $\phi_0\colon \mathcal{R}\rightarrow\mathbb{R}$ defined by
\begin{equation*}
\phi_0(x) = Ae^{\lambda |x|} + B\ \ (1\leq |x| \leq R)
\end{equation*}
where the constants $\lambda < 0$, $A > 0$ and $B$ satisfy the conditions
\begin{equation*}
\left\{\begin{array}{l} Ae^{\lambda} + B = 1\\ Ae^{\lambda R} + B = -1 \end{array}\right.
\end{equation*}
Then for a suitable value of $\lambda < 0$, it holds that
\begin{equation*}
\begin{split}
-\Delta\phi_0 + f(\phi_0) &\leq -\Delta\phi_0 = -A\left(\lambda^2 + \lambda\frac{n-1}{|x|}\right)e^{\lambda |x|} = -\frac{2e^{\lambda |x|}}{e^{\lambda} - e^{\lambda R}}\left(\lambda^2 + \lambda\frac{n-1}{|x|}\right)\\
&\leq -\frac{2e^{\lambda R}}{e^{\lambda} - e^{\lambda R}}\left(\lambda^2 + \lambda\frac{n-1}{|x|}\right) = -\mu < -2\varepsilon
\end{split}
\end{equation*}
in $\mathcal{R}$ for a constant $\mu > 0$, $\phi_0 = -1$ on $\partial B_R$, and $\phi_0 = 1$ on $\partial B_1$, if we take $\delta_0$ such that $0 < \delta_0 = K^{-1}\varepsilon < \frac{1}{2K}\mu$.

Let $\tilde{\phi}$ denote the translation of $\phi_0$ to the ring $B_R(Z)\backslash \bar{B}_1(Z)$. That is $\tilde{\phi}$ satisfies
\begin{equation*}
-\Delta\tilde{\phi} + f(\tilde{\phi}) < -\mu
\end{equation*}
in $B_R(Z)\backslash \bar{B}_1(Z)$ for the constant $\mu > 2\varepsilon$, $\tilde{\phi} = -1$ on $\partial B_R(Z)$, and $\tilde{\phi} = 1$ on $\partial B_1(Z)$.

Now, for each $x$ in $\Omega = B_R(Z)\backslash\bar{B}_1$, define $\tilde{x} = \tau(x)$ in $B_R(Z)\backslash\bar{B}_1(Z)$ in the following way. Write $e = \frac{x}{|x|}$. If
\begin{equation*}
x = (1 - \lambda)e + \lambda q\ \ (0\leq \lambda \leq 1)
\end{equation*}
where $q$ is the point of intersection of the ray from the origin in the direction of $e$ with the sphere $\partial B_R(Z)$, then
\begin{equation*}
\tilde{x} = \tau(x) = (1-\lambda)(Z+e) + \lambda p = Z + (1-\lambda)e + \lambda Re,
\end{equation*}
where $p$ is the point of intersection of the ray from the point $Z$ in the direction of $e$ with the sphere $\partial B_R(Z)$. Clearly, the mapping $x\mapsto \tilde{x}$ is a one-to-one function from $\Omega$ onto $B_R(Z)\backslash\bar{B}_1(Z)$. Suppose $q = te$. Then from $|q - Z| = R$ one can get
\begin{equation*}
t = \sigma(x) := \sqrt{\delta^2\mu^2 + \left(R^2 - \delta^2\right)} - \delta\mu,
\end{equation*}
where $\mu = e\cdot e_1 = x_1/|x|$, and consequently
\begin{equation*}
\lambda = \frac{|x| - 1}{t - 1}.
\end{equation*}
Hence
\begin{equation*}
\tilde{x} = \tau(x) = -\delta e_1 + \left(\frac{t - |x|}{t-1} + \frac{|x| - 1}{t - 1} R\right) e.
\end{equation*}

Finally we define the function $\phi\colon\Omega\rightarrow\mathbb{R}$ by
\begin{equation*}
\phi(x) = \tilde{\phi}(\tilde{x}).
\end{equation*}
We claim that $\phi$ satisfies the conditions
\begin{equation*}
-\Delta\phi + f(\phi) < -\varepsilon
\end{equation*}
in $\Omega$, $\phi = -1$ on $\partial B_R(Z)$, and $\phi = 1$ on $\partial B_1$. In fact, the boundary conditions are obvious. As for the differential inequality, one first writes $\tau = (\tau^1, \tau^2, \hdots, \tau^n)$. Then
\begin{equation*}
\phi_{x_i} = \tilde{\phi}_{\tilde{x}_k}\tau^k_{x_i}
\end{equation*}
and
\begin{equation*}
\phi_{x_ix_j} = \tilde{\phi}_{\tilde{x}_k\tilde{x}_l}\tau^k_{x_i}\tau^l_{x_j} + \tilde{\phi}_{\tilde{x}_k}\tau^k_{x_ix_j}.
\end{equation*}
Here and in the following the summation convention is adopted. Consequently
\begin{equation*}
\phi_{x_ix_i} = \tilde{\phi}_{\tilde{x}_k\tilde{x}_l}\tau^k_{x_i}\tau^l_{x_i} + \tilde{\phi}_{\tilde{x}_k}\tau^k_{x_ix_i}
\end{equation*}
and hence
\begin{equation*}
-\Delta\phi = - <D^2\tilde{\phi}\tau_{x_i}, \tau_{x_i}> - \tilde{\phi}_{\tilde{x}_k}\Delta\tau^k.
\end{equation*}
Decompose $\tau$ as
\begin{equation*}
\tau(x) = x + \psi(x),\ \ \text{where\ }\ \psi(x) = \tau(x) - x.
\end{equation*}
Then
\begin{equation*}
\psi(x) = \tilde{x} - x = -\delta e_1 + \frac{|x|-1}{t-1}\left(R - t\right)e.
\end{equation*}
For any fixed $x\in\Omega$, it is clear that
\begin{equation*}
R - t = \sigma(0) - \sigma(\delta) = -\sigma'(\zeta)\delta
\end{equation*}
for some $\zeta\in (0,\delta)$, and hence
\begin{equation*}
|R - t| \leq 2\delta
\end{equation*}
as $$|\sigma'(\zeta)| = \left|\frac{\delta\mu^2 - \delta}{\sqrt{\delta^2\mu^2 + (R^2 - \delta^2)}} - \mu\right| \leq 2$$ for sufficiently small $\delta$. Moreover, one readily gets
\begin{equation*}
\mu_{x_i} = \frac{\delta_{1i}}{|x|} - \frac{x^1x^i}{|x|^3}
\end{equation*}
and
\begin{equation*}
t_{x_i} = \sigma_{x_i} = \left(\frac{\delta\mu}{\sqrt{\delta^2\mu^2 + (R^2 - \delta^2)}} - 1\right)\delta\mu_{x_i},
\end{equation*}
from which one also gets
\begin{equation*}
\mu_{x_ix_i} = -2\delta_{1i}\frac{x^i}{|x|^3} - \frac{x^1}{|x|^3} + 3\frac{x^1(x^i)^2}{|x|^5}
\end{equation*}
and
\begin{equation*}
t_{x_ix_i} = \frac{R^2 - \delta^2}{\left(\delta^2\mu^2 + (R^2 - \delta^2)\right)^2} \delta^2\mu^2_{x_i} + \left(\frac{\delta\mu}{\sqrt{\delta^2\mu^2 + (R^2 - \delta^2)}} - 1\right)\delta\mu_{x_ix_i}.
\end{equation*}
Clearly, $$\left|\mu_{x_i}\right| \leq \frac{C}{|x|} \leq C\ \ \text{in\ }\ \Omega,$$ and hence $$\left|t_{x_i}\right| \leq C\delta\ \ \text{in\ }\ \Omega.$$

Now
\begin{equation}\label{1deripsi}
\psi_{x_i} = \beta_{x_i}\left(R - t\right)e - \beta t_{x_i}e + \beta\left(R - t\right)e_{x_i},
\end{equation}
where $\beta = \left(|x|-1\right)/\left(t-1\right)$. Evidently $\beta\in [0,1]$ is bounded, and
\begin{equation*}
\left|\beta_{x_i}\right| = \left|\frac{\frac{x_i}{|x|}(t-1) - (|x|-1)t_{x_i}}{(t-1)^2}\right| \leq \frac{1}{t-1} + \frac{|x|-1}{(t-1)^2}C\delta \leq C
\end{equation*}
in $\Omega$. In addition, that
\begin{equation*}
e_{x_i} = \frac{1}{|x|}e_i - \frac{x_i}{|x|^2}e
\end{equation*}
implies $|e_{x_i}|\leq C$ in $|x|\geq 1$. Then one deduces from (\ref{1deripsi}) that
\begin{equation*}
\left|\psi_{x_i}\right| \leq C\delta\ \ \text{in $\Omega$.}
\end{equation*}
Next, one readily gets
\begin{equation}\label{psiii}
\psi_{x_ix_i} = \beta_{x_ix_i}\left(R-t\right)e - \beta t_{x_ix_i}e + \beta\left(R-t\right)e_{x_ix_i} - 2\left(\beta_{x_i}t_{x_i}e + \beta t_{x_i}e_{x_i} - \beta_{x_i}\left(R-t\right)e_{x_i}\right)
\end{equation}
It is clear from the formula of $\mu_{x_ix_i}$ that it is bounded on $\Omega$, which helps to imply from the formula of $t_{x_ix_i}$ that $\left|t_{x_ix_i}\right| \leq C\delta$ on $\Omega$. Meanwhile, one may compute the formula of $\beta_{x_ix_i}$:
\begin{equation*}
\beta_{x_ix_i} = \left(\frac{1}{|x|} - \frac{x^2_i}{|x|^3}\right)\frac{1}{t-1} - \frac{x_i}{|x|}\frac{t_{x_i}}{(t-1)^2} - \frac{x_i}{|x|}\frac{1}{(t-1)^2} + 2\frac{(|x|-1)}{(t-1)^3}t^2_{x_i} - \frac{|x|-1}{(t-1)^2}t_{x_ix_i}.
\end{equation*}
This formula shows that $\left|\beta_{x_ix_i}\right|\leq C$ in $\Omega$ on account of the estimates on $t_{x_i}$ and $t_{x_ix_i}$. Similarly, one gets the formula of $e_{x_ix_i}$
\begin{equation*}
e_{x_ix_i} = -\frac{2x_i}{|x|^3}e_i - \frac{1}{|x|^2}e + 2\frac{x^2_i}{|x|^4}e
\end{equation*}
and deduce from which that $\left|e_{x_ix_i}\right| \leq C$ in $\Omega$. Then the formula (\ref{psiii}) readily implies $\left|\psi_{x_ix_i}\right| \leq C\delta$ on account of the estimates on $R-t$, $\beta$, $e$, $\beta_{x_i}$, $t_{x_i}$, $e_{x_i}$, $\beta_{x_ix_i}$, $t_{x_ix_i}$ and $e_{x_ix_i}$, which in turn implies
$\left|\Delta\psi^k\right| \leq C\delta$ for each $k = 1, \hdots, n$. Computation based on the definition of $\tilde{\phi}$ that
\begin{equation*}
\tilde{\phi}(x) = Ae^{\lambda |x + \delta e_1|} + B
\end{equation*}
and the formulas that determine the values of $A$, $B$ and $\lambda$ helps one to conclude that $\tilde{\phi}_{x_k}$ and $\tilde{\phi}_{x_kx_l}$ are bounded on $\overline{B_R(Z)}\backslash B_1(Z)$. Combining all the preceding estimates, one concludes that
\begin{equation*}
\begin{split}
-\Delta\phi + f(\phi) &= -\Delta\tilde{\phi} + f(\tilde{\phi}) - 2\sum_i<D^2\tilde{\phi}e_i, \psi_{x_i}> - \sum_i<D^2\tilde{\phi}\psi_{x_i}, \psi_{x_i}> - \tilde{\phi}_{x_k}\Delta\psi^k \\
&< -\mu - 2\sum_i<D^2\tilde{\phi}e_i, \psi_{x_i}> - \sum_i<D^2\tilde{\phi}\psi_{x_i}, \psi_{x_i}> - \tilde{\phi}_{x_k}\Delta\psi^k \\
&< -\mu + C\delta \\
&< -\frac{1}{2}\mu, \ \ \text{if we take $K > 2C$}\\
&< -\varepsilon
\end{split}
\end{equation*}
for all $\delta \leq \delta_0$. So the claim is proved.

\subsubsection{Construction of a strict subsolution of $\Delta u=f(u)$ on $\Omega$ satisfying the boundary conditions associated with our problem and the condition $u-\varepsilon \leq v_1 \leq u$ on $\Omega$ }

First replace the subsolution $v_1$ by $\omega_1 :=v_1-C_1\delta$, where $C_1 > \frac{4R}{R-\delta_0}\sup|\nabla u|$. The new function $\omega_1$ satisfies the following conditions:
\begin{equation*}
\left\{\begin{array}{l} u - (\varepsilon + C_1\delta)  \leq \omega_1 \leq u- C_1\delta \ \ \text{in $\Omega$}\\ -\Delta \omega_1 + f(\omega_1) < -(\varepsilon - C_0C_1\delta) < 0\ \ \text{in $\Omega$}\\ \omega_1 < -1\ \ \text{on\ } \partial B_R(Z),\ \ \text{and\ }\ \omega_1 <1\ \ \text{on\ } \partial B_1\end{array}\right.
\end{equation*}
for $C_0 = -\inf_{\mathbb{R}}f'(s) > 0$, if $K$ is sufficiently large.

If one checks carefully our proof in the preceding subsection, it is proved that $-\Delta\tilde{\phi} < -\mu$ and $-\Delta\phi < -\varepsilon$. Then on $\partial\Omega$, $u = \phi$, and in $\Omega$
\begin{equation*}
\Delta(u - \phi) = f(u) - \Delta\phi \leq f(u) - \varepsilon \leq 0.
\end{equation*}
Then the Minimum Principle for super-harmonic functions implies that $u\geq \phi$ on $\bar{\Omega}$.

We are in a position to replace the sub-solution $\omega_1$ by $\tilde{v}_1 := \max\left\{\omega_1, \phi\right\}$ which is also a sub-solution of the problem. Moreover, $\tilde{v}_1$ takes constant values on the exterior and interior spheres respectively. Without any possible confusion, we simply write $v_1$ for $\tilde{v}_1$ in the following. Since $v_1$ differs from $\phi$ on a precompact set, we may mollify it near the boundary of the set. The mollified function $v_1$ verifies $v_1\in C^2(\bar{\Omega})$,
\begin{equation*}
\left\{\begin{array}{l} u - (\varepsilon + 2C_1\delta) \leq v_1 \leq u-\frac{C_1}{2}\delta\ \ \text{in $\Omega$}\\ -\Delta v_1 + f(v_1) < -(\varepsilon - 2C_0C_1\delta) < 0\ \ \text{in $\Omega$}\\ v_1 = -1\ \ \text{on\ } \partial B_R(Z),\ \ \text{and\ }\ v_1 = 1\ \ \text{on\ } \partial B_1\end{array}\right.
\end{equation*}
provided $K$ is sufficiently large.

\subsubsection{Construction of a function $v_{01}$ strict subsolution of $\Delta u=f(u)$ on $\Omega_1$ satisfying the boundary conditions associated with our problem and the condition $u-\varepsilon \leq v_{01} \leq u$ on $\Omega_1$ }

We are ready to define a function $v_0 := v_{01}\in C^2(\overline{\Omega}_1)$ that satisfies
\begin{equation}\label{initial}
\left\{\begin{array}{l} u - C\delta  < v_0 \leq u\ \ \ \text{in $\Omega_1$}\\ -\Delta v_0 + f(v_0) < 0\ \ \text{in $\Omega_1$}\\ v_0 = -1\ \ \text{on $\partial B_{R_1}$ and }\ \ v_0 = 1\ \ \text{on\ }\ \partial B_1\end{array}\right.
\end{equation}
as the initial data for the evolution based on the strict sub-solution$v_1$, where we use and will use in the following $v_0$ for $v_{01}$ to avoid the use of disturbing double subscripts.

For $x\in \overline{\Omega}_1$, if one can write it as
\begin{equation*}
x = (1-\lambda)e + \lambda q,
\end{equation*}
where $e = x/|x|$ and $q = R_1 e$ is the point of intersection of the ray from the origin in the direction of $e$ with the sphere $\partial B_{R_1}$, then one defines
\begin{equation*}
x^* = (1-\lambda)e + \lambda p,
\end{equation*}
where $p$ is the point of intersection of the ray from the origin in the direction of $e$ with the sphere $\partial B_R(Z)$. Clearly, the mapping $x\mapsto x^*$ is a bijection from $\overline{B_{R_1}}\backslash B_1$ onto $\overline{B_R(Z)}\backslash B_1$. Write $p = te$ for $t > 0$. The condition $\left|p + \delta e_1\right| = R$ implies that
\begin{equation*}
t = \sigma(x) := \sqrt{\delta^2\mu^2 + (R^2-\delta^2)} - \delta\mu.
\end{equation*}
Also, we know $\lambda = \frac{|x|-1}{R_1 - 1}$. So
\begin{equation*}
x^* = \varphi(x) := \left(\frac{R_1 - |x|}{R_1 - 1} + \frac{|x| - 1}{R_1 - 1}t\right)e.
\end{equation*}
Set in $\Omega_1$
\begin{equation*}
\psi(x) = \varphi(x) - x = x^* - x = \frac{|x| - 1}{R_1 - 1}\left(t - R_1\right)e.
\end{equation*}
We introduce the notation
\begin{equation*}
\beta(x) = \frac{|x| - 1}{R_1 - 1}.
\end{equation*}
Then
\begin{equation*}
\psi(x) = \beta(x)\left(\sigma(x) - R_1\right)e.
\end{equation*}

Now one can define
\begin{equation*}
v_0(x) = v_1(x^*)\ \ (x \in \overline{\Omega}_1)
\end{equation*}
we claim that $v_0$ satisfies the conditions (\ref{initial}).

The regularity and boundary conditions are evident.

To see that $u - C\delta < v_0 \leq u$ in $\Omega_1$, we write
\begin{equation*}
v_0(x)-u(x)=v_1(x^{*})-u(x)=\big(v_1(x^{*})-u(x^{*})\big)-\big(u(x)-u(x^{*})\big)),
\end{equation*}
which implies
\begin{equation}\label{stab1}
v_0(x) - u(x) \leq -\frac{C_1}{2}\delta + \sup|\nabla u||x-x^*| \leq -\frac{C_1}{2}\delta + \frac{2\delta R}{R_1}\sup|\nabla u| < 0
\end{equation}
and
\begin{equation}\label{stab2}
v_0(x) - u(x) \geq -\varepsilon - 4C_1\delta = -\left(K - 4C_1\right)\delta.
\end{equation}
Here we note that the global gradient estimate of $u$ implies $\sup|\nabla u|$ is controlled by $n$, $R$, and $f$.

Finally, we verify  the differential inequality.

Obviously $\beta$ and $e$ are bounded. The term
\begin{equation*}
\begin{split}
&\ \ \ \ \sigma(x) - R_1 \\
&= \sqrt{\delta^2\mu^2 + \left(R^2 - \delta^2\right)} - \delta\mu - R_1 \\ &= \sqrt{\delta^2\mu^2 + \left(R^2 - \delta^2\right)} - \delta\mu - \left(R - \delta\right) \\ &=: \tau(\delta)
\end{split}
\end{equation*}
for any fixed $x\in \Omega_1$. As $\tau(0) = 0$ and
\begin{equation*}
\left|\tau'(\delta)\right| = \left|\frac{\mu^2\delta - 2\delta}{\sqrt{\delta^2\mu^2 + \left(R^2 - \delta^2\right)}} - \mu + 1\right| \leq C,
\end{equation*}
one concludes
\begin{equation*}
\left|\sigma - R\right| \leq C\delta.
\end{equation*}

One easily gets
\begin{equation*}
\psi_{x_i} = \beta_{x_i}\left(\sigma - R_1\right)e + \beta\sigma_{x_i}e + \beta\left(\sigma - R_1\right)e_{x_i}
\end{equation*}
and
\begin{equation}\label{psiii2}
\psi_{x_ix_i} = \beta_{x_ix_i}\left(\sigma - R_1\right)e + \beta\sigma_{x_ix_i}e + \beta\left(\sigma - R_1\right)e_{x_ix_i} + 2\left(\beta_{x_i}\sigma_{x_i}e + \beta_{x_i}\left(\sigma - R_1\right)e_{x_i} + \beta\sigma_{x_i}e_{x_i}\right).
\end{equation}
Set $\mu(x) = e\cdot e_1 = \frac{x^1}{|x|}$. Then
\begin{equation*}
\mu_{x_i} = \frac{\delta_{1i}}{|x|} - \frac{x^1x^i}{|x|^3}
\end{equation*}
and
\begin{equation*}
\sigma_{x_i} = \frac{\delta^2\mu\mu_{x_i}}{\sqrt{\delta^2\mu^2 + \left(R^2 - \delta^2\right)}} - \delta\mu_{x_i} = \left(\frac{\delta\mu}{\sqrt{\delta^2\mu^2 + \left(R^2 - \delta^2\right)}} - 1\right)\delta\mu_{x_i}.
\end{equation*}
Also
\begin{equation*}
\begin{split}
&\beta_{x_i} = -\frac{1}{R_1-1}\frac{x^i}{|x|}\\
\text{and\ \ } &e_{x_i} = \frac{1}{|x|}e^i - \frac{x^i}{|x|^2}e.
\end{split}
\end{equation*}
As $\left|\mu_{x_i}\right| \leq C$ in $\Omega_1$, it holds $\left|\sigma_{x_i}\right| \leq C\delta$ in $\Omega_1$. Also one observes $\left|\beta_{x_i}\right| \leq C$ and $\left|e_{x_i}\right| \leq C$ in $\Omega_1$. Consequently, it holds
\begin{equation*}
\left|\psi_{x_i}(x)\right| \leq C\delta\ \ \ (x\in\Omega_1).
\end{equation*}

Further computation shows that
\begin{equation*}
\beta_{x_ix_i} = -\frac{1}{(R_1 - 1)|x|} + \frac{x^2_i}{(R_1 - 1)|x|^3}
\end{equation*}
and
\begin{equation*}
e_{x_ix_i} = -\frac{2x^i}{|x|^3}e_i - \frac{1}{|x|^3}x + \frac{3x^2_i}{|x|^5}x,
\end{equation*}
which imply that
\begin{equation*}
\left|\beta_{x_ix_i}\right|,\ \ \left|e_{x_ix_i}\right| \leq C
\end{equation*}
in $\Omega_1$. By computing
\begin{equation*}
\mu_{x_ix_i} = - 2\frac{\delta_{1i}x^i}{|x|^3} - \frac{x^1}{|x|^3} + 3\frac{x^1(x^i)^2}{|x|^5},
\end{equation*}
and
\begin{equation*}
\sigma_{x_ix_i} = \left(\frac{\delta\mu}{\sqrt{\delta^2\mu^2 + \left(R^2 - \delta^2\right)}} - 1\right)\delta\mu_{x_ix_i} + \frac{R^2-\delta^2}{\left( \sqrt{\delta^2\mu^2 + \left(R^2 - \delta^2\right)}\right)^3}\delta^2\mu^2_{x_i},
\end{equation*}
one concludes $\left|\mu_{x_ix_i}\right| \leq C$ in $\Omega_1$ and hence
\begin{equation*}
\left|\sigma_{x_ix_i}(x)\right| \leq C\delta\ \ \ (x\in\Omega_1).
\end{equation*}
The above estimates and the formula (\ref{psiii2}) of $\psi_{x_ix_i}$ imply that
\begin{equation*}
\left|\psi_{x_ix_i}(x)\right| \leq C\delta\ \ \ (x\in\Omega_1)
\end{equation*}

Since
\begin{equation*}
v_{0,x_i} = v_{1,x^*_k}\varphi^k_{x_i}
\end{equation*}
and
\begin{equation*}
v_{0,x_ix_i} = \sum_{k,l}v_{1,x^*_kx^*_l}\varphi^k_{x_i}\varphi^l_{x_i} + \sum_k v_{1,x^*_k} \varphi^k_{x_ix_i},
\end{equation*}
one gets
\begin{equation*}
-\Delta v_0 = -<D^2v_1\varphi_{x_i}, \varphi_{x_i}> - \sum_kv_{1,x^*_k}\Delta\varphi^k.
\end{equation*}
As $\varphi_{x_i} = e_i + \psi_{x_i}$, one further gets from the above formula
\begin{equation*}
-\Delta v_0 = -\Delta v_1 - 2<D^2v_1e_i, \psi_{x_i}> - <D^2v_1\psi_{x_i}, \psi_{x_i}> - \sum_kv_{1,x^*_k}\Delta\psi^k.
\end{equation*}
So
\begin{equation}\label{diff}
\begin{split}
-\Delta v_0 + f(v_0) &= -\Delta v_1 + f(v_1) - 2<D^2v_1e_i, \psi_{x_i}> - <D^2v_1\psi_{x_i}, \psi_{x_i}> - \sum_kv_{1,x^*_k}\Delta\psi^k \\
&< -\left(\varepsilon - 2C_0C_1\delta\right)+C\delta + C\delta \\
&< -C\delta, \ \ \text{for a new constant $C$ if $K$ is sufficiently large.} \\
&< 0
\end{split}
\end{equation}
for all $\delta\leq \delta_0$, on account of the estimates on $e_i$, $\psi_{x_i}$ and $\psi_{x_ix_i}$.

The inequalities in (\ref{stab1}), (\ref{stab2}) and (\ref{diff}) yield to the desired result (\ref{initial}).

\subsubsection{Construction of $w_1(x,t)$ a solution of the parabolic version of our problem on $\Omega_1 \times (0,\infty)$}
Using $v_0$ as the initial data, we are going to solve the following initial-boundary-value problem
\begin{equation}\label{w1est}
\left\{\begin{array}{l}w_t - \Delta w + f(w) = 0\ \text{in\ } \Omega_1\times (0, \infty)\\ w(x,t) = -1\ \ \text{on\ }\partial B_{R_1}\times (0,\infty),\ \ w(x,t) = 1\ \ \text{on\ }\partial B_1\times (0,\infty)\\ w(x,0) = v_0(x)\ \ \text{for\ }x\in\overline{\Omega_1}\end{array}\right.
\end{equation}
For convenience, one sets $\mathcal{D}_1 := \Omega_1\times (0, \infty)$ and let $\partial_p\mathcal{D}_1$ be its parabolic boundary.
\begin{lemma}
There is a solution $w_1$ of the evolution (\ref{w1est}).
\end{lemma}
\begin{pf}
We prove an existence theorem for the following initial-boundary-value problem rewritten from (\ref{w1est}).
\begin{equation}\label{para_exist}
\left\{\begin{array}{ll}w_t - \Delta w + f(w) = 0 &\ \text{in\ }\mathcal{D}_1\\
w(x,t) = v_0(x) &\ \text{on\ } \partial_p\mathcal{D}_1,\end{array}\right.
\end{equation}
where $v_0\in C(\partial_p\mathcal{D}_1)$ is described as before. As $f$ is not proper in the sense it is not a nondecreasing function, one may introduce a function $v(x,t) = e^{-\lambda t}w(x,t)$ in $\mathcal{D}_1$ for a large constant $\lambda >> \frac{2(n+2)}{R^2}$. The function $w$ is a solution of (\ref{para_exist}) if and only if the new function $v$ is a solution of the initial-boundary-value problem
\begin{equation*}
\left\{\begin{array}{ll}v_t - \Delta v + g(t,v) = 0 &\ \text{in\ }\mathcal{D}_1\\
v(x,t) = -e^{-\lambda t}\ \ \text{on\ }\partial B_{R_1}\times (0,\infty),\ \ v(x,t) = e^{-\lambda t}\ \ \text{on\ }\partial B_1\times (0,\infty)\\
v(x,0) = v_0(x) &\ \text{on\ } \bar{\Omega}_1,\end{array}\right.
\end{equation*}
where $g(t,v) = \lambda v + e^{-\lambda t}f(e^{\lambda t}v)$ is a $C^3$ function that is proper, namely $g$ is increasing in $v$. In addition, $g(t,0) = 0$ for any $t$. For simplicity of notation, one may set $\sigma(t)$ be the lateral boundary data of $v$. Writing $w$ for $v$ in the above problem, we are to prove the existence of a solution of the initial-boundary-value problem
\begin{equation}\label{para_exist2}
\left\{\begin{array}{ll}w_t - \Delta w + g(t,w) = 0 &\ \text{in\ }\mathcal{D}_1\\
w(x,t) = \sigma(t) &\ \text{on\ }\left(\partial B_{R_1}\cup\partial B_1\right)\times (0,\infty)\\
w(x,0) = v_0(x) &\ \text{on\ } \bar{\Omega}_1,\end{array}\right.
\end{equation}
The solution of this problem should be well-known. However, as we have not found a proof of the exact problem in the literature, we outline a proof for the reader's convenience. Our proof is different from the usual Perron's method used to attack the existence problem for an elliptic or parabolic equation. Rather, we employed an iterative process to finish the game.

One first picks a function $w^0\in C^2(\bar{\mathcal{D}}_1)$ and proceeds to solve the initial-boundary-value problem
\begin{equation}\label{iteration}
\left\{\begin{array}{ll}w^1_t - \Delta w^1 + g(t,w^0) = 0 &\ \text{in\ }\mathcal{D}_1\\
w^1(x,t) = \sigma(t) &\ \text{on\ }\left(\partial \Omega_1\right)\times (0,\infty)\\
w^1(x,0) = v_0(x) &\ \text{on\ } \bar{\Omega}_1,\end{array}\right.
\end{equation}
for the unknown function $w^1$. This problem can be solved first on the cylinder $\mathcal{D}_{2T} := \Omega_1\times (0, 2T]$ for a small $T$:
\begin{equation*}
\left\{\begin{array}{ll}w^1_t - \Delta w^1 + g(t,w^0) = 0 &\ \text{in\ }\mathcal{D}_{2T}\\
w^1(x,t) = \sigma(t) &\ \text{on\ }\left(\partial \Omega_1\right)\times (0,2T]\\
w^1(x,0) = v_0(x) &\ \text{on\ } \bar{\Omega}_1,\end{array}\right.
\end{equation*}
One then proceeds solving the problem on the cylinder $\Omega_1\times [T,3T]$ with the proper initial-boundary data. The parabolic comparison principle then implies the solutions obtained on the cylinders $\mathcal{D}_{2T}$ and $\Omega_1\times (T,3T]$ coincide on the overlapping part of the two cylinders. And one moves on to the cylinders $\Omega_1\times (2T, 4T]$, $\Omega\times (3T, 5T]$, etc. In the end, one finds a unique solution $w\in C^2(\mathcal{D}_1)$ of (\ref{iteration}) which is $C^2$ up to the vertical boundary. In order to show $w$ is $C^1$ down to the bottom $\Omega_1\times \{t=0\}$, one just differentiates the equation with respect to $t$ to find that $v := w_t$ verifies the conditions
\begin{equation*}
\left\{\begin{array}{ll}v_t - \Delta v + g_t(t,w^0) + g_w(t,w^0)w^0_t = 0 &\ \text{in\ }\mathcal{D}_1\\
v(x,t) = \sigma'(t) &\ \text{on\ }\partial \Omega_1\times (0,\infty)\\
v(x,0) = \Delta v_0(x) - g(0,w^0(x,0)) &\ \text{on\ } \bar{\Omega}_1,\end{array}\right.
\end{equation*}
from which the classical regularity theory of linear equations shows $v$ is continuous down to the bottom. Next, employing the same scheme, one may proceed to solve for each $k = 1, 2, ...$ the initial-boundary-value problem
\begin{equation*}
\left\{\begin{array}{ll}w^{k+1}_t - \Delta w^{k+1} + g(t,w^k) = 0 &\ \text{in\ }\mathcal{D}_1\\
w^{k+1}(x,t) = \sigma(t) &\ \text{on\ }\left(\partial \Omega_1\right)\times (0,\infty)\\
w^{k+1}(x,0) = v_0(x) &\ \text{on\ } \bar{\Omega}_1.\end{array}\right.
\end{equation*}
The functions $w^k$ are $C^2$ up to the lateral sides, and $w_t$ is continuous down to the bottom.

Let $v^k = w^{k+1} - w^k$. Then $v^k$ solves the initial-boundary-value problem
\begin{equation*}
\left\{\begin{array}{ll}v^k_t - \Delta v^k + g(t,w^k) - g(t, w^{k-1}) = 0 &\ \text{in\ }\mathcal{D}_1\\
v^k = 0 &\ \text{on\ } \partial_p\mathcal{D}_1,\end{array}\right.
\end{equation*}
or equivalently,
\begin{equation}\label{eq}
\left\{\begin{array}{ll}v^k_t - \Delta v^k + \tilde{g}(t,x)v^{k-1} = 0 &\ \text{in\ }\mathcal{D}_1\\
v^k = 0 &\ \text{on\ } \partial_p\mathcal{D}_1,\end{array}\right.
\end{equation}
where $\tilde{g}(t,x) = \int^1_0g_w(t, (1-\mu)w^{k-1} + \mu w^k)\,d\mu$.

From here, one easily gets
\begin{equation}\label{bound}
\int_{\Omega_1}\frac{1}{2}\left(v^k(x,T)\right)^2 + \int^T_0\int_{\Omega_1}\left|\nabla v^k\right|^2 = -\int^T_0\int_{\Omega_1}\tilde{g}(t,x)v^kv^{k-1},
\end{equation}
which implies
\begin{equation*}
\frac{1}{2}\int_{\Omega_1}\left(v^k(x,T)\right)^2\,dx \leq \left(\int^T_0\int_{\Omega_1}\frac{1}{2}\left(v^k\right)^2\,dx\,dt\,\int^T_0\int_{\Omega_1} 2\tilde{g}^2\left(v^{k-1}\right)^2\,dx\,dt\right)^{1/2}
\end{equation*}
The latter inequality leads to the estimates
\begin{equation*}
\int^T_0\int_{\Omega_1}\frac{1}{2}\left(v^k\right)^2 \leq CT^2\int^T_0\int_{\Omega_1}\frac{1}{2}\left(v^{k-1}\right)^2,
\end{equation*}
and hence
\begin{equation*}
\int^T_0\int_{\Omega_1}\frac{1}{2}\left(v^k\right)^2 \leq \lambda\int^T_0\int_{\Omega_1}\frac{1}{2}\left(v^{k-1}\right)^2
\end{equation*}
for some $\lambda\in [0, 1)$ if $T$ is small enough. The inequality (\ref{bound}) also gives
\begin{equation*}
\begin{split}
\int^T_0\int_{\Omega_1}\left|\nabla v^k\right|^2 &\leq - \int^T_0\int_{\Omega_1}\tilde{g}(t,x)v^kv^{k-1} \\
&\leq \left(\int^T_0\int_{\Omega_1}\frac{1}{2}\left(v^k\right)^2\,\int^T_0\int_{\Omega_1}\tilde{g}^2\left(v^{k-1}\right)^2\right)^{1/2} \\
&\leq \lambda \int^T_0\int_{\Omega_1}\frac{1}{2}\left(v^{k-1}\right)^2,
\end{split}
\end{equation*}
if one takes the value of $T$ smaller and a new value of $\lambda\in [0, 1)$ if necessary. So $\left\{w^k\right\}$ is a Cauchy sequence with respect to the norm
\begin{equation*}
\|w^k\|_2 = \left(\int^T_0\int_{\Omega_1}\left(w^k\right)^2 + \left|\nabla w^k\right|^2\right)^{1/2}.
\end{equation*}
The equation (\ref{eq}) then implies the boundedness of $w_t$ in the operator norm $\|w_t\|$. As a consequence, a subsequence of $\left\{w^k\right\}$, which we will also denote by $\left\{w^k\right\}$, converges to a certain $w^{\infty}$ in the norm $\|\cdot\|_2$, and the time derivatives $\left\{w^k_t\right\}$ converges weakly to $w^{\infty}_t$. Hence $w^{\infty}$ is a weak solution of (\ref{para_exist2}) on $\bar{\Omega}_1\times [0,T]$. Repeating this process on the time intervals $[\frac{T}{2}, \frac{3T}{2}]$, $[T, 2T]$, $[\frac{3T}{2}, \frac{5T}{2}]$, ..., and employing the parabolic comparison principle, one can find a solution of (\ref{para_exist2}) in $\mathcal{D}_1$. The classical regularity theory then implies $w^{\infty}\in C^2(\mathcal{D}_1)\cap C(\bar{\mathcal{D}}_1)$ (\cite{LSU}, \cite{CLW}, etc). In fact, $w^{\infty}$ is $C^2$ up to the vertical lateral boundary. Moreover, as we did before, one can see $v:=w^{\infty}_t$ solves the linear problem
\begin{equation*}
\left\{\begin{array}{ll}v_t - \Delta v + g_t(t,w^{\infty}) + g_w(t,w^{\infty})v = 0 &\ \text{in\ }\mathcal{D}_1\\
v(x,t) = \sigma'(t) &\ \text{on\ }\partial \Omega_1\times (0,\infty)\\
v(x,0) = \Delta v_0(x) - g(0,w^{\infty}(x,0)) &\ \text{on\ } \bar{\Omega}_1,\end{array}\right.
\end{equation*}
Then $w^{\infty}_t = v$ is continuous down to the bottom. We set $w_1 = e^{\lambda t}w^{\infty}$, and this is the solution we started to obtain. The proof is complete.
\end{pf}
\subsubsection{Convergence of the evolution to a steady state}
We prove the convergence of the evolution (\ref{w1est}) to a steady state.
\begin{lemma}\label{u1}
$$\lim_{t\rightarrow\infty} w_1(x,t) = u_1(x)$$ locally uniformly on $\bar{\Omega}_1$ for some function $u_1$. As a consequence, $u_1$ solves the boundary value problem
\begin{equation*}
\left\{\begin{array}{ll} \lap u = f(u) &\ \text{in $1\leq |x| \leq R_1$}\\ u = 1 &\ \text{on $|x| = 1$}\\ u = -1 &\ \text{on $|x| = R_1$}\end{array}\right.
\end{equation*}
and satisfies
\begin{equation*}
u(x) - C\delta \leq u_1(x) \leq u(x)\ \ \text{in $\Omega_1$.}
\end{equation*}
\end{lemma}
\begin{pf}
Set $z(x,t) = w_{1,t}(x,t)$ on $\bar{\mathcal{D}}_1$. Then $z$ solves the linear initial-boundary-value problem
\begin{equation*}
\left\{\begin{array}{ll}z_t - \Delta z + f'(w_1)z = 0 &\ \text{in\ }\mathcal{D}_1\\
z(x,t) = 0 &\ \text{on\ }\partial \Omega_1\times (0,\infty)\\
z(x,0) = \Delta v_0(x) - f(v_0) &\ \text{on\ } \bar{\Omega}_1,\end{array}\right.
\end{equation*}
Notice that $z\geq 0$ on $\partial_p\mathcal{D}_1$. As $v(x,t)\equiv 0$ is a sub-solution of the above problem with zero initial-boundary data, the parabolic comparison principle implies $z\geq 0$ on $\bar{\mathcal{D}}_1$. Since $u$ is a solution of the evolutionary equation $$u_t - \Delta u + f(u) = 0$$ in $\mathcal{D}_1$ and $u\geq w_1$ on $\partial_p\mathcal{D}_1$, we conclude $w_1(x,t) \leq u(x)$ for all $x\in\bar{\Omega}_1$ and $t \geq 0$. Therefore
\begin{equation*}
\lim_{t\rightarrow +\infty}w_1(x,t) = u_1(x) \leq u(x)
\end{equation*}
monotonically for some function $u_1$ on $\bar{\Omega}_1$. According to either Theorem 3 in \cite{C1} or Theorem 1 in \cite{C2}, it holds that
\begin{equation*}
\|\nabla w_1\|_{L^{\infty}\left(\Omega'\times(0,\infty)\right)} \leq C\left(\|v_0\|_{L^{\infty}(\bar{\Omega})}, \Omega'\right).
\end{equation*}
for any subdomain $\Omega'\subset\subset\Omega_1$. Therefore $w_1(x,t)$ converges to $u_1$ as $t\rightarrow +\infty$ locally uniformly on $\bar{\Omega}_1$. The proof is complete, if one further notices the boundary value of $w_1(x,t)$ is independent of $t$, and the monotonicity of $w_1$ in $t$ along with the fact the initial data $v_0$ satisfies the inequality
$$u(x) - C\delta \leq v_0(x) \leq u(x)$$ in $\Omega_1$.
\end{pf}
\begin{lemma}
$u_1\in C^2(\bar{\Omega}_1)$.
\end{lemma}
\begin{pf}
In the preceding proof, we pointed out that $w_1\in C^2(\bar{\Omega}_1\times(0,\infty))$. As a consequence, $u_1$ is Lipschitz continuous up to the boundary $\partial\Omega_1$. The classical theory of the Possion's equation (e.\,g.\,\cite{GT}) implies $u_1$ is $C^2$ up to the boundary.
\end{pf}
\subsection{Construction of a solution of our problem on the perfect rings $\mathcal{R}$ and $\Omega_2 $ respectively}

Following the same steps we can construct $u_0$ and $u_2$ solutions of our problem on $\mathcal{R}$ and $\Omega_2 $ respectively. we outline the construction of the initial data $v_{00}$ and $v_{02}$.

\begin{enumerate}
\item \label{1} Construct a strict superharmonic function in $\Omega$ satisfying the boundary conditions associated with our problem
\item \label{2} Construct a strict supersolution $v_2$ of $\Delta u=f(u)$ on $\Omega$ satisfying the boundary conditions associated with our problem and the condition $v_2- C\delta \leq u \leq v_2$ on $\Omega$
\item \label {3} Construct a strict supersolution $v_{02}$ of $\Delta u=f(u)$ on $\Omega_2$ satisfying the boundary conditions associated with our problem and the condition  $v_{02}- C\delta \leq u \leq v_{02}$ on $\Omega_2$
\item Extend $v_{02}$ to $\mathcal{R}$ such that $v_{02} \equiv -1$ on $\mathcal{R} \setminus \Omega_1$.
The construction of $v_{00}$ is similar.
\end{enumerate}

The remaining of the argument concerning the existence of a solution of the parabolic problems and the convergence of the evolution, as well as the proof of the above steps, are similar to the ones in the previous subsection. For this reason, we omit the details and just state the results in the following lemmas to avoid making this paper unnecessarily long.

\begin{lemma}
Let  $\mathcal{D}_2 = \Omega_2\times(0,\infty)$. There exists a solution $w_2\in C^2(\bar{\Omega}_2\times(0,\infty))\cap C(\bar{\Omega}_2\times[0,\infty))$ of the initial-boundary-value problem
\begin{equation}
\left\{\begin{array}{ll}w_t - \Delta w + f(w) = 0 &\ \text{in\ }\mathcal{D}_2\\
w(x,t) = v_{02}(x) &\ \text{on\ } \partial_p\mathcal{D}_2,\end{array}\right.
\end{equation}
where $v_{0,2}\in C^2(\bar{\Omega}_2)$ satisfies
\begin{equation}
\left\{\begin{array}{l} u \leq v_{02} \leq u + C\delta  \ \ \text{in $\Omega_2$}\\ -\Delta v_{02} + f(v_{02}) > \varepsilon > 0\ \ \text{in $\Omega_2$}\\ v_{0,2} = -1\ \ \text{on $\partial B_{R_2}$ and }\ \ v_{02} = 1\ \ \text{on\ }\ \partial B_1\end{array}\right.
\end{equation}
\end{lemma}
\begin{lemma}\label{u2}
$$\lim_{t\rightarrow\infty} w_2(x,t) = u_2(x)$$ locally uniformly and monotonically on $\bar{\Omega}_1$. As a consequence, $u_2$ solves the boundary value problem
\begin{equation*}
\left\{\begin{array}{ll} \lap u = f(u) &\ \text{in $1\leq |x| \leq R_2$}\\ u = 1 &\ \text{on $|x| = 1$}\\ u = -1 &\ \text{on $|x| = R_2$}\end{array}\right.
\end{equation*}
and satisfies
\begin{equation*}
u(x) \leq u_2(x) \leq u(x) + C\delta\ \ \text{in $\Omega$.}
\end{equation*}
\end{lemma}
\begin{lemma}
$u_2\in C^2(\bar{\Omega}_2)$.
\end{lemma}

Similarly we have:

\begin{lemma}
Let  $\mathcal{D} = \mathcal{R}\times(0,\infty)$. There exists a solution $w\in C^2(\bar{\mathcal{R}}\times(0,\infty))\cap C(\bar{\mathcal{R}} \times[0,\infty))$ of the initial-boundary-value problem
\begin{equation*}
\left\{\begin{array}{ll}w_t - \Delta w + f(w) = 0 &\ \text{in\ }\mathcal{D}\\
w(x,t) = v_{00}(x) &\ \text{on\ } \partial_p\mathcal{D},\end{array}\right.
\end{equation*}
where $v_{00}\in C^2(\bar{\Omega}_2)$ satisfies
\begin{equation*}
\left\{\begin{array}{l} v_{01} \leq v_{00} \leq v_{02}  \ \ \text{in $\mathcal{R}$}\\ -\Delta v_{00} + f(v_{00}) > \varepsilon > 0\ \ \text{in $\mathcal{R}$}\\ v_{00} = -1\ \ \text{on $\partial \mathcal{R}$ and}\ \ v_{00} = 1\ \ \text{on\ }\ \partial B_1\end{array}\right.
\end{equation*}
\end{lemma}
\begin{lemma}\label{u0}
$$\lim_{t\rightarrow\infty} w(x,t) = u_0(x)$$ locally uniformly and monotonically on $\bar{\Omega}_1$. As a consequence, $u_0$ solves the boundary value problem
\begin{equation*}
\left\{\begin{array}{ll} \lap u = f(u) &\ \text{in $1\leq |x| \leq R$}\\ u = 1 &\ \text{on $|x| = 1$}\\ u = -1 &\ \text{on $|x| = R$}\end{array}\right.
\end{equation*}
and satisfies
\begin{equation*}
u_1(x) \leq u_0(x)\ \ \text{in $\Omega_1$, and\ }\ u_0(x) \leq u_2(x)\ \ \text{in $\mathcal{R}$.}
\end{equation*}
\end{lemma}
\begin{lemma}
$u_0\in C^2(\bar{\mathcal{R}})$.
\end{lemma}
Applying the result of radial symmetry in the preceding section, we conclude that
\begin{theorem}
The solutions $u_i$, $i = 0, 1, 2$, are radially symmetric functions on $\mathcal{R}$ and $\Omega_i$, $i = 1,2$, respectively. In particular, the free boundaries, $\mathcal{F}_i = \partial\left\{u_i>0\right\}$, $i = 0, 1, 2$, are spheres with the center at the origin.
\end{theorem}

\subsection{Comparison and Stability}
The following lemma states the non-degeneracy of $u_2$ in the positive domain.
\begin{lemma}
Let $d(x)$ be the distance from $x$ to $\mathcal{F}_2$. Then
\begin{equation*}
u_2(x) \geq Cd(x)\ \ \ \text{in\ } \left\{u_2 > 0\right\}.
\end{equation*}
\end{lemma}
\begin{pf}
One notices that $u_2$ is super-harmonic in the positive domain $\left\{u_2 > 0\right\}$ and the fact $\mathcal{F}_2$ is a sphere with the origin as its center. Recalling the boundary estimates for a nonnegative harmonic function (e.\,g.\,\cite{C3}, Lemma 6 and proof), one gets the estimate for $u_2$ in the positive domain by comparing $u_2$ to the harmonic function in $\left\{u_2> 0\right\}$ with the same boundary data as $u_2$.
\end{pf}

It is a simple fact that even if two functions are uniformly very close to each other, their boundaries of zero sets, i.\,e.\,the ``\textit{free boundaries}", in general may be far away from each other. Nevertheless, the non-degeneracy of $u_2$ just established helps us to prove in our problem the following lemma that states the free boundary $\mathcal{F}_1$ is indeed close to the other free boundary $\mathcal{F}_2$.
\begin{lemma}\label{trap}
$$dist(\mathcal{F}_1, \mathcal{F}_2) := \sup_{x\in\mathcal{F}_1}dist(x,\mathcal{F}_2) \leq C\delta.$$
\end{lemma}
\begin{pf}
It is known from the previous results, Lemmas \ref{u1} and \ref{u2}, that $u_1\leq u_2 \leq u_1 + C\delta$ on $B_{R_1}\backslash\bar{B}_1$. The non-degeneracy of $u_2$ proved in the preceding lemma implies that
\begin{equation*}
u_2(x) \geq Cd(x)
\end{equation*}
holds on $\mathcal{F}_1$.

On $\mathcal{F}_1$,
\begin{equation*}
u_1(x) +C\delta = C\delta \geq u_2(x) \geq Cd(x),
\end{equation*}
which implies $d(x) \leq C\delta$ for a new constant $C > 0$. That is
\begin{equation*}
dist(\mathcal{F}_1, \mathcal{F}_2) \leq C\delta
\end{equation*}
\end{pf}

We summarize the part of results of the Lemmas \ref{u1}, \ref{u0} and \ref{u2} on the order of the solutions $u_1$, $u_0$, $u$ and $u_2$ on respective domains in the following theorem.
\begin{theorem}
Let $u_i$, $i = 0, 1, 2$, be as constructed in Lemmas \ref{u0}, \ref{u1} and \ref{u2}. 

Then
$u_1\leq u$ in $\Omega_1$, $u\leq u_2$ in $\Omega$, $u_1\leq u_0$ in $\Omega_1$, and $u_0\leq u_2$ in $\mathcal{R}$.

In particular, we have $$|u(x) - u_0(x)| < C\delta\ \ (x\in\Omega\cap\mathcal{R})$$ and the inclusion of the positive sets as stated below.
$$\left\{u_1>0\right\}\subseteq \left\{u > 0\right\} \subseteq \left\{u_2 > 0\right\}\ \ \text{and}$$
$$\left\{u_1>0\right\}\subseteq \left\{u_0 > 0\right\} \subseteq \left\{u_2 > 0\right\}.$$
\end{theorem}
\begin{pf}
The first conclusion is evident from the lemmas mentioned. We need only to point out that $$|u(x) - u_0(x)| < C\delta\ \ (x\in\Omega\cap\mathcal{R})$$ follows from the estimates in the Lemmas \ref{u1} and \ref{u2} and the first conclusion of this theorem. The inclusion of the sets is clear from the first conclusion.
\end{pf}

And in the end by applying Lemma \ref{trap}, we have the desired approximate radial symmetry of $u$.
\begin{theorem}
Let $u$ be as in Theorem \ref{appr_symmetry}, $u_i\ (i = 1, 2)$ be as Lemmas \ref{u1} and \ref{u2}, and $\mathcal{F}$, $\mathcal{F}_i\ (i = 1, 2)$ be their respective free boundaries.

Then
$$dist(\mathcal{F}, \mathcal{F}_0) \leq dist(\mathcal{F}_1, \mathcal{F}_2) < C|Z| = C\delta.$$
\end{theorem}
\begin{pf}
This theorem follows immediately from the inclusion of sets in the preceding theorem and Lemma \ref{trap}.
\end{pf}

The proof of Theorem \ref{appr_symmetry} is complete.

\textbf{Acknowledgement} The authors would like to thank the anonymous referee for many valuable suggestions.

\end{document}